\documentclass{Rsepubli1}
\usepackage{graphicx}
\usepackage{amsmath}
\usepackage{amsfonts}
\usepackage{graphics}
\usepackage{amssymb}
\usepackage{cite}
\usepackage{color}

\usepackage[left=3cm, right=3cm, top=3cm, bottom=3cm
]{geometry}

\newcommand{\e}{^\varepsilon}
\newcommand{\eps}{{\varepsilon}}
\newcommand{\ds}{\displaystyle}

\newcommand{\M}{M\e}

\newcommand{\dist}{\mathrm{dist}}

\renewcommand{\a}{\alpha}
\renewcommand{\b}{\beta}

\newcommand{\cupl}{\bigcup\limits}
\newcommand{\supp}{\mathrm{supp}}
\newcommand{\suml}{\sum\limits}
\newcommand{\intl}{\int\limits}
\newcommand{\liml}{\lim\limits}
\newcommand{\maxl}{\max\limits}
\newcommand{\minl}{\min\limits}

\newcommand{\A}{\mathbf{A}}

\renewcommand{\phi}{\varphi}
\newcommand{\I}{\mathcal{I}(\eps)}
\newcommand{\D}{\mathrm{D}}
\newcommand{\Q}{\mathrm{Q}}

\renewcommand{\l}{\langle}
\renewcommand{\r}{\rangle}

\begin{document}

\title[Homogenization of spectral problem on Riemannian manifold]{Homogenization of spectral problem on Riemannian manifold
consisting of two domains connected by many tubes}
\author[Andrii Khrabustovskyi]{\large{\textbf{Andrii
Khrabustovskyi}}\\\smallskip\normalsize
B.Verkin Institute for Low Temperature Physics and Engineering\\
of the National Academy of Sciences of Ukraine\\\smallskip
\upshape(\texttt{e-mail: andry9@ukr.net})} \maketitle

\begin{abstract}The paper deals with the asymptotic behavior as $\eps\to 0$
of the spectrum of Laplace-Beltrami operator $\Delta\e$ on the
Riemannian manifold $M\e$ ($\mathrm{\dim} M\e=N\geq 2$) depending
on a small parameter $\eps>0$. $M\e$ consists of two perforated
domains which are connected by array of tubes of the length $q\e$.
Each perforated domain is obtained by removing from the fix domain
$\Omega\subset \mathbb{R}^N$ the system of $\eps$-periodically
distributed balls of the radius $d\e=\bar{o}(\eps)$. We obtain a
variety of homogenized spectral problems in $\Omega$, their type
depends on some relations between $\eps$, $d\e$ and $q\e$. In
particular if the limits $\liml_{\eps\to 0}q\e$ and
$\liml_{\eps\to 0}\ds{(d\e)^{N-1}q\e \eps^{-N}}$ are positive then
the homogenized spectral problem contains the spectral parameter
in a nonlinear manner, and its spectrum has a sequence of
accumulation points.
\end{abstract}

\parindent=1.5pc
\section*{Introduction}
The aim of this paper is to study the asymptotic behavior as
$\eps\to 0$ of the spectrum of Laplace-Beltrami operator
$\Delta\e$ (with Dirichlet boundary conditions) on the
$N$-dimensional Riemannian manifold $M\e$ $(N\geq 2)$ depending on
small parameter $\eps>0$. The manifold $M\e$ is embedded in
$\mathbb{R}^{N+1}$. It is constructed in the following way. Let
$\Omega$ be a bounded smooth domain in $\mathbb{R}^N$ and let
$\left\{D_i\e\right\}_{i}$ be a system of disjoint balls ("holes")
of the radius $d\e$ distributed $\eps$-periodically in $\Omega$.
Denote $\Omega\e=\Omega\setminus\cupl_i D_i\e$. Then $M\e$
consists of two parallel perforated sets $\Omega_1\e$ and
$\Omega_2\e$ (each of them is the copy of $\Omega\e$) and the set
$\left\{G_i\e\right\}_{i}$ of cylinders of the length $q\e$ and
the radius $d\e$ (the cylinder $G_i\e$ connects the boundaries of
$i$-th holes in $\Omega_1\e$ and $\Omega_2\e$):
$$M\e=\Omega_1\e\cup\left(\cupl_{i}G_i\e\right)\cup\Omega_2\e.$$
The manifold $M\e$ is presented on Fig.1. We equip $M\e$ by
Riemannian metric $g\e$ which is induced on $M\e$ by Euclidean
metric in $\mathbb{R}^{N+1}$. More precise description of $M\e$
will be specified later in Section \ref{sec1}.
\begin{figure}[h]
\begin{center}
\scalebox{0.3}[0.3]{\includegraphics{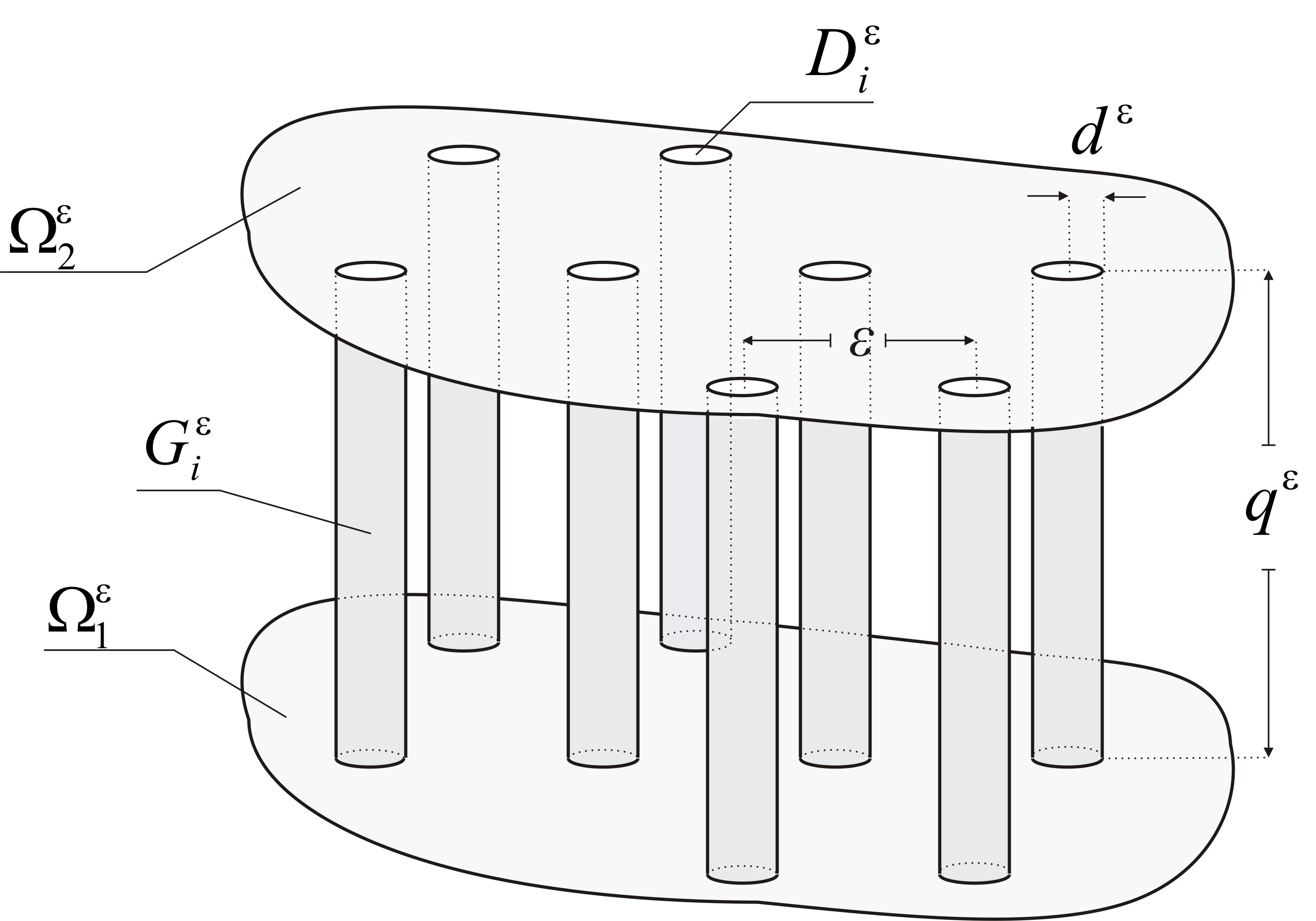}}
 \caption{The manifold $M\e$}
\end{center}
\end{figure}

We denote by $\sigma(-\Delta\e)=\{\lambda\e_{m}\}_{m\in
\mathbb{N}}$
 the sequence of eigenvalues of $-\Delta\e$
(here they are renumbered in the increasing order and repeated
according to their multiplicity). By $\left\{u\e_m\right\}_{m\in
\mathbb{N}}$ we denote the system of corresponding eigenfunctions
that are chosen orthonormal in ${L_2(M\e)}$.

Our goal is to find the homogenized spectral problem in $\Omega$
whose spectrum is a limit of $\sigma(-\Delta\e)$ as $\eps\to 0$.
Let us note that we choose the Dirichlet boundary conditions only
for the sake of definiteness, all results are still valid for
Neumann or mixed boundary conditions too.

In one special case of the relationship between $\eps$, $d\e$ and
$q\e$ this problem was studied in \cite{DalMaso} (see the remark
after Theorem \ref{th4} below). In the current work we impose much
more weaker restrictions on $d\e$ and $q\e$ comparing with the
work \cite{DalMaso}.

Firstly the homogenization problem on Riemannian manifolds of
complex microstructure was studied  in \cite{Khrus1}. The
investigations in \cite{Khrus1} are motivated by the problem to
describe asymptotic behavior of colored particles moving in the
domain with small obstacles when the number of obstacles tends to
infinity: it turns out that this problem can be reduced to the
homogenization of diffusion equation on some Riemannian manifold
depending on small parameter $\eps$.

The next works in this direction were devoted to the
homogenization of semi-linear parabolic equations and their
attractors \cite{Khrus2}, the homogenization of harmonic vector
fields \cite{Khrus3} and the homogenization of Maxwell equations
\cite{Khrus4}. The works \cite{Khrus3,Khrus4} are related to the
general relativity (according to Wheeler \cite{Wheeler} such
manifolds can be interpreted as models of the Universe). Some
applications of the homogenization theory on manifolds were also
presented in \cite{Khrab3}.

The asymptotic behavior of the spectrum of Laplace-Beltrami
operator on Riemannian manifolds of complex microstructure firstly
was studied in \cite{Notar}. In this work the manifold $M\e$
consists of some fixed manifold (possibly without a boundary) and
an increasing number of attached thin handles. Close problems were
also considered in \cite{DalMaso,Khrab1}. The same problem on the
manifolds with \textit{one} attached handle of small thickness was
studied in \cite{Anne,Chavel1,Chavel2} (see also the survey
\cite{Kuwae}, where the convergence of spectra is studied on
various Riemannian manifolds depending on a small parameter but
the dependence on this parameter has essentially another nature
comparing with homogenization problems). The spectral problems on
manifolds consisting of a fixed manifold with increasing number of
attached small spherical manifolds ("bubbles") were studied in
\cite{Khrab4,Khrab2}.

In the works
\cite{DalMaso,Khrus1,Khrus2,Khrus3,Khrus4,Khrab1,Khrab4,Notar} it
is assumed that the radiuses $d\e$ of holes are of order
$\eps^{N\over N-2}$ if $N>2$ or $\exp\left(-a/\eps^2\right)$ if
$N=2$ (incidentally, the homogenization of Dirichlet BVP for
Poisson equation in the perforated domain with such holes leads to
the appearance of the potential like term in the homogenized
equation, see e.g. \cite{Chio,March}). Also in the works
\cite{DalMaso,Khrab1,Notar} it is supposed that the length of the
attached tubes tends to zero as $\eps\to 0$.

As it was mentioned above in the current work we impose much more
weaker restrictions on the sizes of holes and tubes: we suppose
that $d\e=\bar{o}(\eps)$ as $\eps\to 0$, and the total volume of
tubes and the length $q\e$ of tubes are bounded uniformly in
$\eps$ ($\eps<\eps_0$). For more precise statement see the
conditions (\ref{main_cond}) below. For example if
$d\e=\mathbf{d}\eps^{\alpha}$, $q\e=\mathbf{q}\eps^\beta$
($\mathbf{d},\mathbf{q}$ are positive constants) then these
conditions are valid iff $\alpha>1$, $\alpha(N-1)+\beta-N\geq 0$
and $\beta\geq 0$.

Under these assumptions we obtain a variety of qualitatively
different types of the homogenized spectral problem. It turns out
that the type of homogenized spectral problem depends essentially
on the limits (\ref{main_cond}), (\ref{D}), (\ref{Q}) below. The
most attention is devoted to the case when the both limits
$q=\liml_{\eps\to 0} q\e$ and $p=\liml_{\eps\to
0}\ds{(d\e)^{N-1}q\e\eps^{-N}}$ exist and are positive. In this
case the spectrum $\sigma(-\Delta\e)$ converges to the set
$\mathcal{A}=\sigma\left(\mathrm{A}(\lambda)\right)\cup\left(\cupl_{n\in
\mathbb{N}}\left\{(\pi n)^2q^{-2}\right\}\right)$, where
$\sigma\left(\mathrm{A}(\lambda)\right)$ is the spectrum of some
operator pencil $\mathrm{A}(\lambda)$: each operator
$\mathrm{A}(\lambda)$ acts in $[L_2(\Omega)]^2$ and contains the
spectral parameter $\lambda$ in a non-linear manner (see Theorem
\ref{th1} below). The spectrum of $\mathrm{A}(\lambda)$ consists
of the sequence $\left\{\lambda_{m}^n\right\}_{m,n\in \mathbb{N}}$
of isolated eigenvalues with finite multiplicity
 such that for
fixed $n$ the subsequence $\{\lambda_{m}^n\}_{m\in \mathbb{N}}$
belongs to the open segment $\left((\pi (n-1))^2q^{-2}, (\pi
n)^2q^{-2}\right)$ and $\lambda_{m}^n\underset{m\to\infty}\nearrow
(\pi n)^2q^{-2}$.

If $p=0$ (but $q$ is still positive) the pencil
$\mathrm{A}(\lambda)$ becomes a linear (see Theorem \ref{th2}).

In the case $q=0$ the spectrum $\sigma(-\Delta\e)$ converges to
the spectrum of some homogenized operator $\mathrm{A}$ acting
either in $L_2(\Omega)$ or in $[L_2(\Omega)]^2$ and having purely
discrete spectrum (see Theorems \ref{th3}-\ref{th4}).
\medskip

\noindent\textit{Remark.} Such a structure of the spectrum of
homogenized problem as in the case $q>0$, $p>0$ is also
characteristic for the problems posed on co-called thick
junctions. Thick junctions are domains with highly oscillating
boundary: they consist of a junction body and a great number of
attached thin domains located along a joining zone on the surface
of the junction body. Boundary-value problems in thick junctions
were studied by many authors (see, e.g.,
\cite{Blan1,Blan2,Melnyk1,Melnyk2,Melnyk3,Khrus5}). In particular
in the work \cite{Melnyk2} the asymptotic behavior as $\eps\to 0$
of eigenvalues and eigenfunctions of the Neumann problem is
investigated on the junction $\Omega\e\subset \mathbb{R}^2$
consisting of two domains connected by an $\eps$-periodic system
of thin strips of fixed length. Just as in the current work the
spectrum of the homogenized problem in \cite{Melnyk2} consists of
the sequence of isolated eigenvalues with finite multiplicity and
of the points $\{P_n\}_{n\in \mathbb{N}}$ that divide the
eigenvalues into countably many subsequences convergent to the
corresponding point $P_n$.

Another problem that leads to such structure of the spectrum is
consider in \cite{Zhikov}. Here the author investigates the
asymptotic behavior as $\eps\to 0$ of the spectrum of operator
$\mathrm{A}\e=\mathrm{div}\left(a\e(x)\nabla\right)$ in fixed
bounded domain with coefficient $a\e(x)$ that degenerates as
$\eps\to 0$ on some disperse periodic set. The operator
$\mathrm{A}\e$ corresponds to double-porosity media, at present
there is a great number of works related to this field (see, e.g.,
the books \cite{Chechkin,March} and references therein).\medskip

The outline of the paper is the following. In Section \ref{sec1}
we describe precisely the structure of the manifold $\M$ and
formulate the main results of the paper (Theorems
\ref{th1}-\ref{th4}) describing the Hausdorff convergence of
$\sigma(-\Delta\e)$ as $\eps\to 0$. Also we illustrate these
results on the example mentioned above (i.e.
$d\e=\mathbf{d}\eps^{\alpha}$, $q\e=\mathbf{q}\eps^\beta$). In
Section \ref{sec2} we present some auxiliary technical lemmas
which are used in the proof of main results. In Section \ref{sec3}
we prove Theorems \ref{th1}-\ref{th4}. The proof is based on the
substitution of suitable test functions into the variational
formulation of the spectral problem as in the energy method using
for classical homogenization problems (see e.g. the books
\cite{ChioDon,Tartar,Sanch}). Finally in Section \ref{sec4} we
present the results on a number-by-number convergence of the
eigenvalues, i.e. convergence as $\eps\to 0$ of $\lambda_m\e$ for
fixed number $m$ (Theorems \ref{th6}, \ref{th8}, \ref{th10}), and
the convergence of eigenfunctions $u\e_m$ (Theorems \ref{th7},
\ref{th9}).

\section{\label{sec1}Setting of problem and main results}

Let $\Omega$ be a bounded smooth domain in $\mathbb{R}^N$ $(N\geq
2)$ and let $D_i\e\ (i\in \I\subset \mathbb{Z}^N)$ be a system of
disjoint balls ("holes") of the radius $d\e$ with centers at
$x_i\e=i\cdot\eps\ (i\in\mathbb{Z}^N)$
 such that $D_i\e\subset\Omega$ and $\dist(x_i\e,\partial\Omega)\geq {\eps\over 2}$. Here $\I$ stands for
corresponding set of multiindexes $i$. We denote
$$\Omega\e=\Omega\setminus\left(\cupl_{i\in\I} D_i\e\right).$$
In $\mathbb{R}^{N+1}$ we consider the following sets (below $x\in
\mathbb{R}^N,\ z\in \mathbb{R},\ (x,z)\in \mathbb{R}^{N+1}$):
\begin{gather*}
\Omega_1\e=\big\{(x,z)\in \mathbb{R}^{N+1}:\ x\in\Omega\e,\
z=0\big\},\quad \Omega_2\e=\big\{(x,z)\in \mathbb{R}^{N+1}:\
x\in\Omega\e,\ z=q\e\big\},\\ G_i\e=\big\{(x,z)\in
\mathbb{R}^{N+1}:\ x\in\partial D_i\e,z\in[0,q\e]\big\},
\end{gather*}
where $q\e$ is a positive number.

Finally we obtain the set $M\e$ in $\mathbb{R}^{N+1}$ consisting
of two perforated domains $\Omega_1\e$ and $\Omega_2\e$ which are
connected by the set of cylinders $G_i\e$:
$$M\e=\Omega_1\e\cup\left(\cupl_{i\in \I}G_i\e\right)\cup\Omega_2\e.$$
We denote by $\tilde x$ points of $M\e$. Also we denote
$$S_{1i}\e=\big\{(x,z)\in \mathbb{R}^{N+1}: x\in\partial D_i\e,\ z=0\big\},\quad
S_{2i}\e=\big\{(x,z)\in \mathbb{R}^{N+1}: x\in\partial D_i\e,\
z=q\e\big\}.$$

Clearly $\M$ can be covered by a system of charts and suitable
local coordinates $\{x_1,\dots,x_N\}\mapsto \tilde x\in\M$ can be
introduced. In particular in a small neighbourhood
 of $S_{1i}\e$ we introduce them as follows. Let
$(\phi_1,\dots,\phi_{N-1},r)$ be the spherical coordinates in
$\Omega_1\e$ with the origin at $x_i\e$. Here
$\phi_1,\dots,\phi_{N-1}$ are the angular coordinates
($\phi_1\in[0,2\pi],\phi_{j}\in[0,\pi]\ (j={2,\dots,N-1})$), $r$
($r\geq d\e$) is the distance to $x_i\e$ (that is $r=d\e$ for the
points of $S_{1i}\e$). Let $(\phi_1,\dots,\phi_{N-1},z)$ be the
cylindrical coordinates in $G_i\e$. We set $x_j\e=\phi_j\
(j={1,\dots,N-1})$, $x_N=r-d\e\ (x_N\geq 0)$ for
${\tilde{x}}\in\Omega_1\e$ and $x_N=-z$ ($x_N\leq 0$) for
$\tilde{x}\in G_i\e$. Similarly local coordinates can be introduce
in a small neighbourhood of $S_{2i}\e$.

Therefore we obtain the $N$-dimensional differential manifold
$M\e$. If the point $\tilde x$ belongs to $\Omega_k\e$ ($k=1,2$)
we assign to $\tilde x$ a pair $(x,k)$, where $x$ is a
corresponding point in $\Omega\e$. If the point $\tilde x$ belongs
to $G_i\e$ ($i\in\I$) we assign to $\tilde x$ a pair $(\phi,z)$,
where $\phi=(\phi_1,\dots,\phi_{N-1})$ are the angular
coordinates, $z\in[0,q\e]$. The boundary of $M\e$ consists of the
exterior boundaries of $\Omega_1\e$ and $\Omega_2\e$, i.e.
$\partial M\e=\cupl_{k=1,2}\big\{\tilde x=(x,k)\in \Omega_k\e:
x\in\partial\Omega\}.$

The Euclidean metrics in $\mathbb{R}^{N+1}$ induces on the
manifold $M\e$ the Riemannian metrics
$g\e=\left\{g\e_{\a\b}\right\}_{\a,\b=\overline{1,N}}$. It is
clear that the metrics $g\e$ is continuous and piecewise-smooth
(it is smooth everywhere outside the $(N-1)$-dimensional spheres
$S_{ki}\e$ ($k=1,2,\ i\in\I$)). In a small neighbourhood of
$S_{1i}\e$ ($i\in \I$) the components of $g\e$ have the following
form in the local coordinates $(x_1,\dots,x_N)$ introduced above:
\begin{gather*}
g_{\a\b}\e=\delta_{\a\b}\cdot\begin{cases}
(x_N+d\e)^2\prod\limits_{j=\a+1}^{N-1}\sin^{2}x_j,&x_N\geq 0,\\
(d\e)^2\prod\limits_{j=\a+1}^{N-1}\sin^{2}x_j,&x_N<0,
\end{cases}\ \a=\overline{1,N-1},\qquad
g_{n\beta}=\delta_{n\beta}
\end{gather*}
(for $\a=N-1$ we set $\prod\limits_{j=\a+1}^{N-1}\sin^2 x_j$:=1).
Here $\delta_{\a\b}$ is the Kronecker's delta.

Let us introduce the following functional spaces:
\begin{itemize}
\item $L_2(M\e)$ be the Hilbert space of square integrable (with
respect to the volume measure) functions on $M\e$. The scalar
product and norm are defined by
$$(u,v)_{L_2(M\e)}={\intl_{M\e} u\bar{v}d\tilde
x},\quad  \|u\|_{L_{2}(M\e)}=\sqrt{(u,u)_{L_{2}(M\e)}},$$ where
$d\tilde x=\sqrt{\det g\e}d x_1{\dots}d x_N$ is the volume measure
on $M\e$;

\item $H^1(M\e)$ be the Hilbert space of square integrable
functions on $M\e$ with gradient from $L_2(M\e)$. The scalar
product and norm are defined by
$$(u,v)_{H^1(M\e)}={\intl_{M\e}\bigg(\nabla\e u\cdot\nabla\e \bar{v}+ u\bar{v}\bigg)
d\tilde x},\quad  \|u\|_{H^{1}(M\e)}=\sqrt{(u,u)_{H^1(M\e)}},$$
where $\nabla\e u\cdot\nabla\e \bar{v}$ is the scalar product of
the vector fields $\nabla\e u$ and $\nabla\e \bar{v}$ with respect
to the metrics $g\e$. In local coordinates $\nabla\e
u\cdot\nabla\e \overline{v}=\ds\suml_{\a,\b=1}^N
g_\eps^{\a\b}\ds{\partial u\over\partial x_\a}{\partial
\bar{v}\over\partial x_\b}$, where $g_\eps^{\a\b}$ are the
components of the tensor inverse to $g^\eps_{\alpha\beta}$;

\item ${H_0^1}(M\e)$ be the subspace of $H^1(M\e)$ consisting of
functions $u$: $u|_{\partial M\e}=0$.
\end{itemize}

It is well-known (see e.g. \cite{Taylor}) that for any $f\e\in
L_2(M\e)$ there exists the unique $u_f\e\in H_0^1(M\e)$ such that
\begin{gather*}
(\nabla\e u_f\e,\nabla\e
v\e)_{L_2(M\e)}=(f\e,v\e)_{L_2(M\e)},\quad \forall v\e\in
H_0^1(M\e).
\end{gather*}
Thus we have the operator $\mathrm{T}\e$ that acts in $L_2(M\e)$
and is defined by the formula $\mathrm{T}\e f\e=u_f\e$. This
operator is compact and self-adjoint. We denote
$\Delta\e=-(\mathrm{T}\e)^{-1}$. The operator $\Delta\e$ is called
\textit{Laplace-Beltrami operator (with Dirichlet boundary
conditions)}. In local coordinates it has the following form:
$$\Delta\e={1\over \sqrt{\det{ g\e}}}\suml_{\a,\b=1}^N {\partial\over
\partial x_\a}\left(\sqrt{\det g\e} g_\eps^{\a\b}{\partial\over\partial
x_\b}\right)
$$

$\Delta\e$ is the self-adjoint operator with purely discrete
spectrum. We denote by $\sigma(-\Delta\e)=\{\lambda\e_{m}\}_{m\in
\mathbb{N}}$ the sequence of eigenvalues of $-\Delta\e$ written in
the increasing order and repeated according to their multiplicity:
$$0<\lambda_1\e\leq\lambda_2\e\leq{\dots}\leq\lambda_m\e\leq{\dots}\underset{m\to 0}\to\infty.$$
By $\left\{u\e_m\right\}_{m\in \mathbb{N}}$ we denote the system
of corresponding eigenfunctions such that
$(u_\a\e,u_\b\e)_{L_2(M\e)}=\delta_{\a\b}$.

Our goal is to describe the asymptotic behavior of
$\sigma(-\Delta\e)$ as $\eps\to 0$. As it was mentioned in the
introduction we choose the Dirichlet boundary conditions only for
the sake of definiteness.
\medskip

\noindent\textit{Remark.} We have noted above that $g\e$ is the
piecewise-smooth metrics. Nevertheless $g\e$ can be easily
approximated by the smooth metrics $g^{\eps\delta}$ that differs
from $g\e$ only in small $\delta(\eps)$-neighborhoods of
$S\e_{ki}$ while in this neighborhoods $g^{\eps}$ and
$g^{\eps\delta}$ are sufficiently close (see e.g. \cite{Khrus1}
for the exact construction). Let $\Delta^{\eps\delta}$ be
Laplace-Beltrami operator corresponding to the metrics
$g^{\eps\delta}$. If $\delta(\eps)$ converges to $0$ sufficiently
fast as $\eps\to 0$ then the limits as $\eps\to 0$ of the
spectrums $\sigma(-\Delta^{\eps\delta})$ and
$\sigma(-\Delta^{\eps})$  are the same. This can be proved using
for example the double-sided inequality in the end of section
\textit{\textit{"Outils"}} in \cite{Anne}. However it is more
convenient to carry out the proof of the results for the
piecewise-smooth metrics $g\e$.\medskip

We will solve our problem under the following assumptions:
\begin{gather}
\label{main_cond} \liml_{\eps\to 0} \ds{d\e\over\eps}=0,\quad
\exists\liml_{\eps\to 0}{(d\e)^{N-1} q\e\over\eps^N}=p,\quad
\exists\liml_{\eps\to 0}q\e=q,\quad p,q\in [0,\infty).
\end{gather}
Obviously the condition $p<\infty$ implies that the total volume
of cylinders $G_i\e$ is bounded uniformly in $\eps$
($\eps<\eps_0$).

In the simplest situation $d\e=\mathbf{d}\eps^{\alpha}$,
$q\e=\mathbf{q}\eps^\beta$ ($\mathbf{d},\mathbf{q}>0$ are
constants) conditions (\ref{main_cond}) are valid iff $\alpha>1$,
$\alpha(N-1)+\beta-N\geq 0$ and $\beta\geq 0$. This example (for
$N>2$) will be discussed after Theorems \ref{th1}-\ref{th4} below.

In order to describe the behavior as $\eps\to 0$ of
$\sigma(-\Delta\e)$ we use a concept of Hausdorff convergence.

\smallskip \noindent\textbf{Definition.}\textit{ Let $\mathcal{A}\e\subset
\mathbb{R}$ be the set depending on the positive parameter $\eps$.
We say that $\mathcal{A}\e$ converges as $\eps\to 0$ in the
Hausdorff sense to the set $\mathcal{A}_0$ if the following
conditions hold:
\begin{gather}\tag{A}\label{ah}
\text{if }\lambda\e\in\mathcal{A}\e\text{ and }\liml_{\eps\to
0}\lambda\e=\lambda_0\text{ then }\lambda_0\in
\mathcal{A}_0,\\\tag{B}\label{bh} \text{for any }\lambda_0\in
\mathcal{A}_0\text{ there is }\lambda\e\in\mathcal{A}\e\text{ such
that }\liml_{\eps\to 0}\lambda\e=\lambda_0.
\end{gather}
}

Now we are able to formulate the main results of the paper.
Starting from the case $q>0$, we introduce the following operator
pencil $\mathrm{A}(\lambda)$, $\lambda\in
\mathbb{\mathbb{C}}\setminus\cupl_{n\in \mathbb{N}}\left\{(\pi
n)^2q^{-2}\right\}$: the operator $\mathrm{A}(\lambda)$ acts in
$[L_2(\Omega)]^2$, it is defined by the operation
\begin{gather*}
\mathrm{A}(\lambda)=\left(\begin{matrix}-\Delta
+\ds{p\omega\sqrt\lambda\over q\tan (q\sqrt{\lambda})}&
-\ds{p\omega\sqrt\lambda\over q\sin (q\sqrt{\lambda})}\\
-\ds{p\omega\sqrt\lambda\over q\sin (q\sqrt{\lambda})}&-\Delta
+\ds{p\omega\sqrt\lambda\over q\tan
(q\sqrt{\lambda})}\end{matrix}\right)-\lambda \mathrm{I}
\end{gather*}
and by  the definitional domain
$\mathcal{D}(\mathrm{A}(\lambda))=\left\{\mathrm{U}\in
[H^2(\Omega)]^2,\ \mathrm{U}|_{\partial\Omega}=0\right\}$. Here by
$\omega$ we denote the volume of $(N-1)$-dimensional unit sphere,
$\mathrm{I}$ is the identical operator. By
$\sigma\left(\mathrm{A}(\lambda)\right)$ we denote the
\textit{spectrum} of the pencil $\mathrm{A}(\lambda)$, i.e. the
set of such $\widehat{\lambda}\in
\mathbb{\mathbb{C}}\setminus\cupl_{n\in \mathbb{N}}\left\{(\pi
n)^2q^{-2}\right\}$ that the operator
$\mathrm{A}(\widehat{\lambda})$ does not have a bounded inverse
operator. A number $\widehat{\lambda}$ is called an
\textit{eigenvalue} of the operator pencil $\mathrm{A}(\lambda)$
if $\mathrm{A}(\widehat{\lambda})\mathrm{U}=0$ for some
$\mathrm{U}\not= 0$.

\begin{theorem}{}\label{th1} Suppose that $q>0$ and $p>0$. Then as
$\eps\to 0$ the spectrum $\sigma(-\Delta\e)$ converges in the
Hausdorff sense to the set
$\mathcal{A}=\sigma\left(\mathrm{A}(\lambda)\right)\cup\left(\cupl_{n\in
\mathbb{N}}\left\{(\pi n)^2q^{-2}\right\}\right)$.

The spectrum $\sigma(\mathrm{A}(\lambda))$ consists of the
isolated eigenvalues $\lambda_{m}^n$ ($m,n\in \mathbb{N}$) with
finite multiplicity which are distributed on the positive semiaxis
in the following way:
\begin{gather}\label{struc1}
\forall n\in \mathbb{N}:\quad(\pi (n-1))^2
q^{-2}<\lambda_{1}^n\leq\lambda_{2}^n\leq{\dots}\leq\lambda_{m}^n\leq{\dots}\underset{m\to\infty}\to
(\pi n)^2q^{-2}
\end{gather}
\end{theorem}

In the last section we will perfect the result of Theorem
\ref{th1} proving that $\forall m\in \mathbb{N}$\ \
$\lambda_m\e\underset{\eps\to 0}\to\lambda_m^1$ (Theorem
\ref{th8}).

\begin{theorem}{}\label{th2} Suppose that $q>0$ and $p=0$. Then as
$\eps\to 0$ the spectrum $\sigma(-\Delta\e)$ converges in the
Hausdorff sense to the set
$\mathcal{A}=\sigma\left(\mathrm{A}\right)\cup\left(\cupl_{n\in
\mathbb{N}}\left\{(\pi n)^2q^{-2}\right\}\right)$. Here
$\sigma\left(\mathrm{A}\right)$ is the spectrum of the operator
$\mathrm{A}$ that acts in $[L_2(\Omega)]^2$ and is defined by the
operation
\begin{gather}\label{A3}
\mathrm{A}=-\left(\begin{matrix}\Delta & 0\\0&\Delta
\end{matrix}\right)
\end{gather}
and by  the definitional domain
$\mathcal{D}(\mathrm{A})=\left\{\mathrm{U}\in [H^2(\Omega)]^2,\
\mathrm{U}|_{\partial\Omega}=0\right\}$.
\end{theorem}

In the last section we will perfect the result of Theorem
\ref{th2} proving that if the number $m$ is sufficiently large
then $\lambda_m\e\underset{\eps\to 0}\to\pi^2 q^{-2}$ (Theorem
\ref{th10}).
\medskip

In the case $q=0$ we additionally suppose that the following
limits exist:
\begin{gather}
\label{D} r=\liml_{\eps\to 0}{(d\e)^{N-1} \over \eps^N q\e},\quad
\D=\liml_{\eps\to 0}\ds{\mathrm{D}\e\over \eps^{N}},\quad
r,\D\in[0,\infty]
\end{gather}
where $$\mathrm{D}\e=\begin{cases} |\ln
d\e|^{-1},&N=2,\\(d\e)^{N-2},& N>2.\end{cases}$$ If $\D\in
(0,\infty)$ we suppose that the following limit exists:
\begin{gather}\label{Q}
\mathrm{Q}=\liml_{\eps\to 0}\begin{cases}\ds{q\e\over d\e|\ln
d\e|},&N=2,\\ \ds{q\e\over d\e},&N>2,
\end{cases}\quad \Q\in [0,\infty]
\end{gather}

\begin{theorem}{}\label{th3} Suppose that $q=0,\ r=\infty,\
\mathrm{D}=\infty$. Then as $\eps\to 0$ the spectrum
$\sigma(-\Delta\e)$ converges in the Hausdorff sense to the
spectrum $\sigma(\mathrm{A})$ of operator $\mathrm{A}$ that acts
in $L_2(\Omega)$ and is defined by the operation
\begin{gather}\label{A1}
\mathrm{A}=-\left(1+{1\over 2}{p\omega}\right)^{-1}\Delta
\end{gather}
and by the definitional domain
$\mathcal{D}(\mathrm{A})=\left\{u\in H^2(\Omega),\
u|_{\partial\Omega}=0\right\}$.
\end{theorem}

\begin{theorem}{}\label{th4} Suppose that $q=0$ and either
$r<\infty,\ \D=\infty$ or $\D<\infty$. Then as $\eps\to 0$ the
spectrum $\sigma(-\Delta\e)$ converges in the Hausdorff sense to
the spectrum of operator $\mathrm{A}$ that acts in
$[L_2(\Omega)]^2$ and is defined by the operation
\begin{gather}\label{A2}
\mathrm{A}=\left(\begin{matrix}-\Delta +\mathrm{V}&
-\mathrm{V}\\-\mathrm{V}&-\Delta +\mathrm{V}\end{matrix}\right)
\end{gather}
and by  the definitional domain
$\mathcal{D}(\mathrm{A})=\left\{\mathrm{U}\in [H^2(\Omega)]^2,\
\mathrm{U}|_{\partial\Omega}=0\right\}$. Here $\mathrm{V}$ is the
operator of multiplication by constant
\begin{gather}\label{VVV}
\mathrm{V}=\begin{cases}r\omega,&0<r<\infty,\ \D=\infty,\\
\ds{(N-2)\omega\D \over 2+(N-2)\Q},&0<\D<\infty,\ \Q<\infty,\
N>2,\\\ds {2\pi\D\over 2+\Q},&0<\D<\infty,\ \Q<\infty,\
N=2,\\0,&(r=0,\ \D=\infty)\ \text{or}\ (0<\D<\infty,\ \Q=\infty)\
\text{or}\ (\D=0).
\end{cases}
\end{gather}
\end{theorem}
\medskip

\noindent\textit{Remark.} In the case $N>2$, $\liml_{\eps\to
0}\ds{d\e\over \eps}=0$ and $\ds{q\e\over d\e}=\mathrm{const}$
(i.e. $G_i\e$ is the $d\e$-homothetic image of the fixed cylinder)
our problem was also investigated in \cite[Chapter 2]{DalMaso} by
using $\Gamma$-convergence technique. In this case $d\e$ and $q\e$
satisfy the conditions of Theorem \ref{th3} (if $\D=\infty$) or
the conditions of Theorem \ref{th4} (if $\D<\infty$). The results
obtained in the current work agree with the results obtained in
\cite{DalMaso}.
\medskip

\noindent\textbf{Example.} We consider the example mentioned
above: let $d\e=\mathbf{d}\eps^{\alpha}$,
$q\e=\mathbf{q}\eps^\beta$ ($\mathbf{d},\mathbf{q}>0$ are
constants). Also let $N>2$.

Let us consider the coordinate plane $(\a,\b)$ (see Fig.2). On
this plane we mark some important points: $A=(1,\infty),\
B=(1,1),\ C=\left({N\over N-1},0\right),\ D=(\infty,0),\
E=\left({N\over N-2},{N\over N-2}\right)$, $F=\left({N\over
N-2},\infty\right),\ G=\left({N\over N-2},0\right)$.
\begin{figure}[h]
\begin{center}
\scalebox{0.4}[0.4]{
\includegraphics{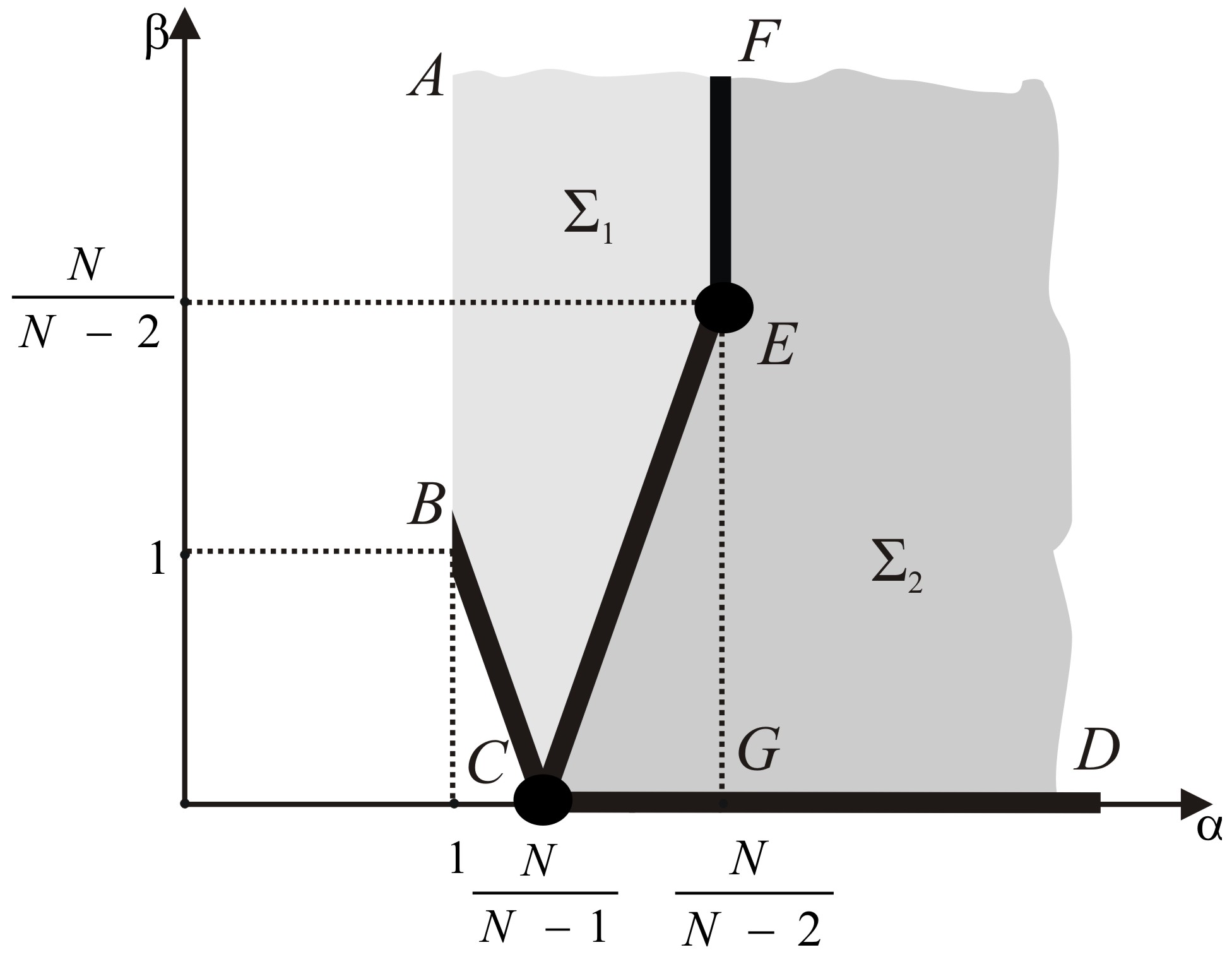}} \caption{Plain $(\alpha,\beta)$}
\end{center}
\end{figure}

Since conditions (\ref{main_cond}) hold iff $\alpha>1$,
$\alpha(N-1)+\beta-N\geq 0$ and $\beta\geq 0$ then we are
restricted by the set involving the open unbounded domain whose
boundary is the polyline $ABCD$, the open segment $(B,C)$ and the
ray $[C,D)$. It is easy to see that Theorems \ref{th1}-\ref{th4}
describe the behavior of $\sigma(-\Delta\e)$ for all
$(\alpha,\beta)$ from this set. Indeed:
\begin{itemize}
\item In the point $C$ we have $q=\mathbf{q},\
p=\mathbf{d}^{N-1}\mathbf{q}$. This case is described by Theorem
\ref{th1}.

\item On the open ray $(C,D)$ we have $q=\mathbf{q},\ p=0$. This
case is described by Theorem \ref{th2}.

\item In the open domain $\Sigma_1$ whose boundary is the polyline
$ABCEF$ we have $q=0$, $r=\infty$, $\D=\infty$. This case is
described by Theorem \ref{th3}. Here $p=0$ and therefore the
homogenized operator is $-\Delta$.

\item On the open segment $(B,C)$ we have $q=0,\ r=\infty,\
\D=\infty$. This case is also described by Theorem \ref{th3} but
here $p=\mathbf{d}^{N-1}\mathbf{q}>0$ and the homogenized operator
is $-\left(1+{1\over
2}\mathbf{d}^{N-1}\mathbf{q}\omega\right)^{-1}\Delta$.

\item In the open segment $(C,E)$ we have either $q=0,\
r=\mathbf{d}^{N-1}\mathbf{q}^{-1},\ \D=\infty$, on the ray $[E,F)$
we have $q=0,\ \D=\mathbf{d}^{N-2},\ \Q<\infty$. These two cases
are described by Theorem \ref{th4} (by the way in the point $E$ we
have $\Q=\mathbf{d}^{-1}\mathbf{q}$, on the open ray $(E,F)$ we
have $\Q=0$). In this case $\mathrm{V}>0$.

\item On the open domain $\Sigma_2$ whose boundary is the polyline
$FECD$ we have $r=0,\ \D=\infty$ (within the triangle $CEG$) or
$\D=0$ (in the open domain whose boundary is the polyline $FGD$)
or $\D=\mathbf{d}^{N-2},\ \Q=\infty$ (on the open segment
$(E,G)$). This case is described by Theorem \ref{th4} but
here $\mathrm{V}=0$, the homogenized operator is\ $-\left(\begin{matrix}\Delta &0\\
0&\Delta\end{matrix}\right)$.
\end{itemize}

We remark that the spectrums of homogenized operators in
$\Sigma_1$ and $\Sigma_2$ coincide but the multiplicity of each
eigenvalue in $\Sigma_2$ is of two times grater then its
multiplicity in $\Sigma_1$. This difference will be taken into
account in Section~\ref{sec4} where a number-by-number convergence
of the eigenvalues is studied.

Also we remark that in the case $N=2$ we have $\D=\infty$ for any
$\alpha$. Therefore in this case in order to cover all types of
homogenized problems we also have to consider the radius $d\e$
that tends to zero faster then $\eps^\alpha$, $\forall\alpha$. For
example if $d\e=\exp(-a/\eps^2)$ ($a\in (0,\infty)$) then
$\D=a^{-1}$.

\section{\label{sec2}Auxiliary results}
In this section we obtain some technical lemmas which are used in
the proof of Theorems \ref{th1}-\ref{th4}.

Let us introduce the list of notations. Recall that if the point
$\tilde x$ belongs to $\Omega_k\e$ ($k=1,2$) we assign to $\tilde
x$ a pair $(x,k)$, where $x$ is a corresponding point in
$\Omega\e$;  if the point $\tilde x$ belongs to $G_i\e$ ($i\in\I$)
we assign to $\tilde x$ a pair $(\phi,z)$, where
$\phi=(\phi_1,{\dots},\phi_{N-1})$ are the angular coordinates,
$z\in[0,q\e]$.

\begin{itemize}
\item $\square_{i}\e$ be the cube in $\mathbb{R}^N$ with the
center at $x_i\e$, side-length $\eps$ and edges which are parallel
to the coordinate axes; \item $B_{ki}\e=\left\{\tilde x=(x,k)\in
\Omega_k\e:\ d\e\leq |x-x_i\e|\leq\eps/2\right\}$; \item
$C_{ki}\e=\left\{\tilde x=(x,k)\in \Omega_k\e:\
|x-x_i\e|=\eps/2\right\}$; \item $S_{i}\e[\tau]=\left\{\tilde
x=(\phi,z)\in G_i\e:\ z=\tau\right\},\ \tau\in[0,q\e]$\quad (that
is $S_{1i}\e=S_{i}\e[0],\ S_{2i}\e=S_{i}\e[q\e]$); \item $\l
u\r_B$\ (where either $B\subset M\e$ or $B\subset\Omega$) be the
mean value of the function $u\in L_2(B)$. That is $\ds\l
u\r_B={1\over |B|}\intl_{B}ud\tilde x$, where $|B|$ is the volume
of $B$; \item $\mathrm{S}_{N-1}$ be the $(N-1)$-dimensional unit
sphere,\ $d\phi=\left(\prod\limits_{k=1}^{N-1}
\sin^{k-1}\phi_{k}\right)d\phi_1{\dots}d\phi_{N-1}$ be the volume
measure on $\mathrm{S}_{N-1}$; \item $C,\ C_1,\ C_2,$ etc. be
generic positive constants independent of $\eps$.
\end{itemize}
We introduce the operators $\widehat{X}\e:C^1(\Omega)\to
L_2(\Omega)$, $\widehat{\square}\e:L_2(\Omega)\to L_2(\Omega)$,
$\widehat{B}_k\e:L_2(M\e)\to L_2(\Omega)$,
$\widehat{S}_k\e:L_2(M\e)\to L_2(\Omega)$ by the following
formulae
\begin{gather}\label{oper1}
\begin{matrix}
[{\widehat{X}}\e u](x)=\begin{cases}u(x_i\e),&
x\in\square_i\e\\0,&\Omega\setminus\cupl_{i\in \I}\square_i\e;
\end{cases}&
[{\widehat{\square}}\e u](x)=\begin{cases}\l u\r_{\square_{i}\e},&
x\in\square_i\e\\0,&\Omega\setminus\cupl_{i\in \I}\square_i\e;
\end{cases}\\\ \\
[\widehat{B}_k\e u\e](x)=\begin{cases}\l u\e\r_{B_{ki}\e},&
x\in\square_i\e\\0,&\Omega\setminus\cupl_{i\in \I}\square_i\e;
\end{cases}&
[\widehat{S}_k\e u\e](x)=\begin{cases}\l u\e\r_{S_{ki}\e},&
x\in\square_i\e\\0,&\Omega\setminus\cupl_{i\in \I}\square_i\e.
\end{cases}
\end{matrix}
\end{gather}
Also we introduce an extension operators $\Pi_k\e:H^1(M\e)\to
H^1(\Omega)$ ($k=1,2$) such that
\begin{gather*}
\text{if }x\in \Omega\e\text{ then } [\Pi_k\e u\e](x)=u\e(\tilde
x),\ \tilde x=(x,k)\in\Omega_k\e;\\ \forall u\e\in
H^1(\Omega_k\e):\ \|\Pi_k\e u\e\|_{H^1(\Omega)}\leq
C\|u\e\|_{H^1(\Omega_k\e)}.
\end{gather*}
It is well-known (see e.g. \cite{ACDP,CSJP}) that such operators
exist.

\begin{lemma}\label{lm1}
Let $u\e\in H^1(M\e)$. Then the following inequalities hold:
\begin{gather}
\label{in_lm1_1}\text{{\bf\rm I}.\qquad} \left|\l
u\e\r_{S_{ki}\e}-\l u\e\r_{B_{ki}\e}\right|^2\leq C{\|\nabla
u\e\|^2_{L_2(B_{ki}\e)}\over \D\e},\quad k=1,2,\ i\in\I,\\
\label{in_lm1_3}\text{{\bf\rm II}.\qquad} \left|\l \Pi_k\e
u\e\r_{\square_{i}\e}-\l u\e\r_{B_{ki}\e}\right|^2\leq
C\eps^{2-N}\|\nabla\Pi_k\e u\e\|^2_{L_2(\square_i\e)},\quad
k=1,2,\ i\in\I,\\ \label{in_lm1_2}\text{{\bf\rm III}.\qquad}
\left|\l u\e\r_{S_i\e[\tau_1]}-\l
u\e\r_{S_i\e[\tau_2]}\right|^2\leq \omega^{-1}\|\nabla\e
u\e\|^2_{L_2(G_i\e)}{q\e\over (d\e)^{N-1}},\quad
\forall\tau_1,\tau_2\in [0,q\e],\ i\in\I.
\end{gather}
\end{lemma}

\noindent\textbf{Proof.} I.\ Let us fix $k$ and $i$, and let us
introduce a spherical coordinates $(\phi,r)$ in $B_{ki}\e$. Here
$r$ is a distance to $x_i\e$ ($r\geq d\e$),
$\phi=(\phi_1,{\dots},\phi_{N-1})$ are the angular coordinates.
Let $x=(\phi,d\e)\in S_{ki}\e$, $y=(\phi,r)\in B_{ki}\e$. We have
\begin{gather*}
u\e(x)-u\e(y)=-\intl_{d\e}^{r}{\partial
u\e(\phi,\tau)\over\partial\tau}d\tau.
\end{gather*}
Then we multiply this equality by $r^{N-1}d\phi dr$, integrate
from $d\e$ to $\eps/2$ (with respect to $r$) and over
$\mathrm{S}_{N-1}$ (with respect to $\phi$), divide by
$|B_{ki}\e|$ and square. Using Cauchy inequality we obtain:
\begin{gather*}
\left|\l u\e\r_{S_{ki}\e}-\l
u\e\r_{B_{ki}\e}\right|^2={1\over|B_i\e|^2}\left|\intl_{d\e}^{\eps/2}r^{N-1}
\intl_{\mathrm{S}_{N-1}}\intl_{d\e}^{r}{\partial
u\e(\phi,\tau)\over\partial\tau} (\phi,\tau)d\tau d\phi
dr\right|^2\leq\\\leq
{\omega\over|B_i\e|^2}\left(\ds{\eps\over2}-d\e\right)\ds\intl_{d\e}^{\eps/2}r^{2(N-1)}dr
\intl_{\mathrm{S}_{N-1}}\intl_{d\e}^{\eps/2}\left|{\partial
u\e(\phi,\tau)\over\partial\tau}\right|^2\tau^{N-1}d\tau
d\phi\intl_{d\e}^{\eps/2}{d\tau\over\tau^{N-1}}\leq C{\|\nabla
u\e\|^2_{L_2(B_{ki}\e)}\over \D\e}.
\end{gather*}

II.\ The inequality (\ref{in_lm1_3}) is a particular case of Lemma
2.1 from \cite{Khrab2}.

III.\ Let $\tilde x=(\phi,\tau_1)\in S_i\e[\tau_1]$, $\tilde
y=(\phi,\tau_2)\in S_i\e[\tau_2]$. Then $ \ds u\e(\tilde
x)-u\e(\tilde y)=-\intl_{\tau_1}^{\tau_2}\ds{\partial
u\e(\phi,\tau)\over\partial\tau}d\tau,$ and we obtain
\begin{multline*}
\left|\l u\e\r_{S_i\e[\tau_1]}-\l
u\e\r_{S_i\e[\tau_2]}\right|^2\leq\\\leq
\omega^{-1}|\tau_1-\tau_2|\intl_{\mathrm{S}_{N-1}}\intl_{\tau_1}^{\tau_2}
\left({\partial u\e(\phi,\tau)\over\partial\tau}\right)^2
 d\tau d\phi_1{\dots}d\phi_{N-1}\leq \omega^{-1}\|\nabla u\e\|^2_{L_2(G_i\e)}{q\e\over
(d\e)^{N-1}}.\quad \square
\end{multline*}

\begin{corollary}\label{cor1}
Let $u\e\in H^1(M\e)$, $\|\nabla\e u\e\|^2_{L_2(M\e)}<C$ and
$\Pi_k\e u\e \underset{\eps\to 0}\to u_k$ strongly in
$L_2(\Omega)$ ($k=1,2$). Then we have
\begin{eqnarray}\label{in_cor1_2}
\text{if\ }\D=\infty:\quad &\widehat{S}_k\e u\e\underset{\eps\to
0}\to u_k,\\\label{in_cor1_1} & \widehat{B}_k\e
u\e\underset{\eps\to 0}\to u_k,\\\label{in_cor1_3}\text{\ if\
}\ds\liml_{\eps\to 0}{(d\e)^{N-1}\over \eps^N}=0:\quad &
\ds{(d\e)^{N-1}\over\eps^N}\|\widehat{S}_k\e
u\e\|^2_{L_2(\Omega)}\underset{\eps\to 0}\to 0,\\\label{in_cor1_4}
\text{if\ }r=\D=\infty:\quad& u_1=u_2\text{\ \ {\rm and}\
}\liml_{\eps\to
0}\ds\suml_{i\in\I}\|u\e\|^2_{L_2(G_i\e)}=p\omega\|u_1\|^2_{L_2(\Omega)}.
\end{eqnarray}
\end{corollary}

\noindent\textbf{Proof.} We present the proof only for the
statement (\ref{in_cor1_2}) (another statements are proved
similarly using Lemma \ref{lm1}). One has:
\begin{multline}\label{4terms}
\|\widehat{S}_k\e u\e-u_k\e\|_{L_2(\Omega)}\leq
\left(\|\widehat{S}_k\e u\e-\widehat{B}_k\e
u\e\|_{L_2(\Omega)}+\|\widehat{B}_k\e u\e-\widehat{\square}\e
\Pi_k\e
u\e\|_{L_2(\Omega)}+\right.\\\left.+\|\widehat{\square}\e\Pi_k\e
u\e-\Pi_k\e u\e\|_{L_2(\Omega)}+\|\Pi_k\e
u\e-u_k\|_{L_2(\Omega)}\right).
\end{multline}
Due to the inequality (\ref{in_lm1_1}) the first term in
(\ref{4terms}) tends to zero if $\D=\infty$:
\begin{gather*}
\|\widehat{S}_k\e u\e-\widehat{B}_k\e
u\e\|^2_{L_2(\Omega)}=\eps^N\suml_{i\in\I}\left|\l
u\e\r_{S_{ki}\e}-\l u\e\r_{B_{ki}\e}\right|^2\leq C{\eps^N\over
\D\e}\|\nabla\e u\e\|^2_{L_2(\Omega_k\e)}\underset{\eps\to 0}\to
0.
\end{gather*}
In a similar way inequality (\ref{in_lm1_3}) implies that the
second term also tends to zero. The third term tends to zero by
virtue of the Poincare inequality for the cube $\square_i\e$. And
finally the last term tends to zero by the given data. Thus the
statement (\ref{in_cor1_2}) is proved.\quad $\square$

\begin{lemma}\label{lm2}  Let $q=p=0$. Let $u\e\in H^1(M\e)$, $\|\nabla\e u\e\|^2_{L_2(M\e)}<C$.  Then
$$\liml_{\eps\to 0}\suml_{i\in\I}\|u\e\|^2_{L_2(G_i\e)}=0.$$

\end{lemma}

\noindent\textbf{Proof} follows directly from the following
inequality: for $\forall u\e\in H^1(M\e)$
\begin{gather}
\label{in_q0} \|u\e\|^2_{G_i\e}\leq C\left\{(q\e)^2\|\nabla
u\e\|^2_{L_2(G_i\e)}+{{q\e}(d\e)^{N-1}\over\eps^N}\|u\e\|^2_{L_2(B_i\e)}+
{q\e(d\e)^{N-1}\over \D\e}\|\nabla
u\e\|^2_{L_2(B_i\e)}\right\}.\end{gather} This inequality is
proved in \cite{Khrab1} (Lemma 2.2) for $N=2$. For $N>2$ the proof
is fully similar.\ $\square$

\begin{lemma}\label{lm3} Let $q>0$. Let $\lambda\e\in\sigma(-\Delta\e)$, $u\e$
be the corresponding eigenfunction such that
$\|u\e\|_{L_2(M\e)}=1$. Suppose that $\lambda\e\underset{\eps\to
0}\to\lambda_0\notin \cupl_{n\in \mathbb{N}}\left\{(\pi n)^2
q^{-2}\right\}$ and $\Pi\e_k u\e\underset{\eps\to 0}\to u_k\in
H_0^1(\Omega)$ ($k=1,2$) strongly in $L_2(\Omega)$.

Then if $p>0$ one has
\begin{gather}\label{in_lm3_1}
\liml_{\eps\to 0}\suml_{i\in\I}\|u\e\|^2_{L_2(G_i\e)}=
p\omega\bigg\{k_1\left(\|u_1\|_{L_2(\Omega)}^2+
\|u_2\|_{L_2(\Omega)}^2\right)+2k_2(u_1,u_2)_{L_2(\Omega)}\bigg\},
\end{gather}
where $k_1=\ds{q\sqrt{\lambda_0}-\sin (q\sqrt{\lambda_0})\cos
(q\sqrt{\lambda_0})\over
2q\sqrt{\lambda_0}\sin^2(q\sqrt{\lambda_0})}$,
$k_2=-\ds{q\sqrt{\lambda_0}\cos(q\sqrt{\lambda_0})-\sin
(q\sqrt{\lambda_0})\over
2q\sqrt{\lambda_0}\sin^2(q\sqrt{\lambda_0})}$.

If $p=0$ one has
\begin{gather}\label{in_lm3_2}
\liml_{\eps\to 0}\suml_{i\in\I}\|u\e\|^2_{L_2(G_i\e)}=0.
\end{gather}

\end{lemma}

\textbf{Proof.} We introduce on $G_i\e$ the function
$v_i\e(\phi,z)=\l u\e\r_{S_i\e[z]}$ (it is clear that $v_i\e$ is
independent of $\phi$). By the Poincare inequality
\begin{gather}\label{puan}
\suml_{i\in\I}\|u\e-v_i\e\|^2_{G_i\e}\leq C(d\e)^2\suml_{i\in\I}
\|\nabla\e u\e\|^2_{L_2(G_i\e)}\leq C(d\e)^2\lambda\e.
\end{gather}

Since $-\Delta\e u\e=\lambda\e u\e$ then it is easy to see that
\begin{gather}
\label{bvp1} -(v_i\e)''=\lambda\e v_i\e,\ z\in (0,q\e),\quad
v_i\e(0)=\l u\e\r_{S_{1i}\e},\quad v_i\e(q\e)=\l u\e\r_{S_{2i}\e}.
\end{gather}

So long as $\lambda\e\underset{\eps\to 0}\to\lambda_0\notin
\cupl_{n\in \mathbb{N}}\left\{(\pi n)^2 q^{-2}\right\}$ and
$q\e\underset{\eps\to 0}\to q$ then for sufficiently small $\eps$
$\lambda\e\not\in \cupl_{n\in \mathbb{N}}\left\{{(\pi n)^2
(q\e)^{-2}}\right\}$. Therefore the problem (\ref{bvp1}) has the
unique solution
\begin{gather*}
v_i\e(z)=A_i\e\sin(z\sqrt{\lambda\e})+B_i\e\cos(z\sqrt{\lambda\e}
),\text{ where  } A_i\e={\l u\e\r_{S_{2i}\e}-\l
u\e\r_{S_{1i}\e}\cos(q\e\sqrt{\lambda\e})\over
\sin(q\e\sqrt{\lambda\e})},\ B_i\e=\l u\e\r_{S_{1i}\e}.
\end{gather*}

Direct computations show that
\begin{gather*}
\|v_i\e\|^2_{L_2(G_i\e)}=\omega(d\e)^{N-1}q\e\bigg\{k\e_1\left(\l
u\e\r^2_{S_{1i}\e}+\l u\e\r^2_{S_{2i}\e}\right)+2k\e_2\l
u\e\r_{S_{1i}\e}\cdot\l u\e\r_{S_{2i}\e} \bigg\},
\end{gather*}
where $\ds k_1\e=\ds{q\sqrt{\lambda\e}-\sin
(q\e\sqrt{\lambda\e})\cos(q\e\sqrt{\lambda\e})\over
2q\e\sqrt{\lambda\e}\sin^2(q\e\sqrt{\lambda\e})}$,
$k_2=-\ds{q\sqrt{\lambda\e}\cos(q\e\sqrt{\lambda\e})-\sin
(q\e\sqrt{\lambda\e})\over
2q\e\sqrt{\lambda\e}\sin^2(q\e\sqrt{\lambda\e})}$. Therefore
\begin{multline}
\label{in_lm3_last}
\suml_{i\in\I}\|v_i\e\|^2_{L_2(G_i\e)}=\\=\omega{(d\e)^{N-1}q\e\over\eps^N}\bigg\{k\e_1
\left(\|\widehat{S}_1\e u\e\|^2_{L_2(\Omega)}+\|\widehat{S}_2\e
u\e\|^2_{L_2(\Omega)}\right)+2k\e_2(\widehat{S}_1\e
u\e,\widehat{S}_2\e u\e)_{L_2(\Omega)}\bigg\}.
\end{multline}
Then (\ref{in_lm3_1}) follows directly from (\ref{in_lm3_last}),
(\ref{puan}) and Corollary \ref{cor1} (see (\ref{in_cor1_2})).
Similarly (\ref{in_lm3_2}) follows from (\ref{in_lm3_last}),
(\ref{puan}) and (\ref{in_cor1_3}). Lemma is proved.\qquad
$\square$

\section{\label{sec3}Proof of the main Theorems}
\subsection{Proof of Theorem \ref{th1}}
\textbf{Step 1.} Firstly we prove that condition (\ref{ah}) of the
Hausdorff convergence holds. Let $\lambda\e\in\sigma(-\Delta\e)$
and $\lambda\e\underset{\eps\to 0}\to\lambda_0$. If $\lambda_0\in
\cupl_{n\in \mathbb{N}}\left\{{(\pi n)^2 q^{-2}}\right\}$ then
(\ref{ah}) is proved. Therefore we are interested in the case
$\lambda_0\notin \cupl_{n\in \mathbb{N}}\left\{{(\pi n)^2
q^{-2}}\right\}$.

Let $u\e$ be the eigenfunction that correspond to $\lambda\e$ and
$\|u\e\|_{L_2(M\e)}=1$ (and therefore $\|\nabla\e
u\e\|^2_{L_2(M\e)}=\lambda\e$). Since the functions $u\e$ are
bounded  in $H_0^1(M\e)$ uniformly in $\eps$ then $\Pi_k\e u\e$
($k=1,2$) are also bounded  in $H_0^1(M\e)$ uniformly in $\eps$.
Therefore due to the embedding theorem there exists a subsequence
(still denoted by $\eps$) such that
\begin{gather*}
\Pi_k\e u\e\underset{\eps\to 0}\to u_k\in H_0^1(\Omega)\quad
(k=1,2)\quad \text{\ strongly in }L_2(\Omega)\text{\ and weakly in
}H^1(\Omega).
\end{gather*}
By Lemma \ref{lm3} we have
$$1=\liml_{\eps\to 0}\|u\e\|^2_{L_2(M\e)}=\|u_1\|^2_{L_2(\Omega)}+
\|u_2\|^2_{L_2(\Omega)}+p\omega
\bigg\{k_1\left(\|u_1\|^2_{L_2(\Omega)}+
\|u_2\|^2_{L_2(\Omega)}\right)+2k_2(u_1,u_2)_{L_2(\Omega)}\bigg\}$$
and therefore
$\mathrm{U}=\left(\begin{matrix}u_1\\u_2\end{matrix}\right)\not=
0$. We prove that $\mathrm{A}(\lambda_0)\mathrm{U}=0$.

For an arbitrary $w\e\in H^1_0(M\e)$ we have:
\begin{gather}\label{int_ineq}
\intl_{M\e}\nabla\e u\e\cdot\nabla\e w\e d\tilde
x-\lambda\e\intl_{M\e}u\e w\e d\tilde x=0.
\end{gather}
Let us introduce the following test function $w\e$:
\begin{gather}\label{we}
w\e(\tilde x)=\begin{cases}w_k(x)+\ds\suml_{i\in
\I}\bigg(w_k(x_i\e)-w_k(x)\bigg)\cdot\phi\left({|x-x_i\e|\over
d\e}\right),&\tilde x=(x,k)\in \Omega_k\e,\ k=1,2,\\
v_i\e(z),&\tilde x=(\phi,z)\in G_i\e,\ i\in\I.
\end{cases}
\end{gather}
Here $w_k(x)\in C_0^\infty(\Omega)$ ($k=1,2$) are arbitrary
functions, $\phi(r):[0,\infty)\to \mathbb{R}$ is a smooth positive
function equal to $1$ as $r\leq 1$ and equal to $0$ as $r\geq 2$,
$v_i\e(z)$ is defined by the formula:
\begin{gather*}
v_i\e(z)=A_i\e\sin(z\sqrt{\lambda\e})+B_i\e\cos(z\sqrt{\lambda\e}
),\text{ where }
A_i\e={w_2(x_i\e)-w_1(x_i\e)\cos(q\e\sqrt{\lambda\e})\over
\sin(q\e\sqrt{\lambda\e})},\ B_i\e=w_1(x_i\e)
\end{gather*}
(we suppose that $\eps$ is sufficiently small so that
$\lambda\e\not\in \cupl_{n\in \mathbb{N}}\left\{{(\pi n)^2
(q\e)^{-2}}\right\}$). It is easy to see that
\begin{gather*}
 -(v_i\e)''=\lambda\e v_i\e,\ z\in (0,q\e),\quad
v_i\e(0)=w_1(x_i\e),\quad v_i\e(q\e)=w_2(x_i\e).
\end{gather*}

We denote $\phi_i\e=\phi\left({|x-x_i\e|/d\e}\right)$.

Substituting this $w\e$ into (\ref{int_ineq}) and taking into
account that
\begin{gather*}
\intl_{G_i\e}\bigg(\nabla\e u\e\cdot\nabla\e w\e-\lambda\e u\e
w\e\bigg)d\tilde x=(d\e)^{N-1}\omega\left(-\l
u\e\r_{S_{1i}\e}{\partial v_i\e(0)\over\partial z}+\l
u\e\r_{S_{2i}\e}{\partial v_i\e(q\e)\over\partial z}\right),
\end{gather*}
we obtain
\begin{multline}\label{plugg1}
0=\suml_{k=1,2}\intl_{\Omega\e}\bigg(\nabla (\Pi_k\e
u\e)\cdot\nabla w_k -\lambda\e (\Pi_k\e u\e)
w_k\bigg)dx+\\+\suml_{i\in
\I}{(d\e)^{N-1}\omega\sqrt{\lambda\e}\over\sin(q\e\sqrt{\lambda\e})}\bigg(\l
u\e\r_{S_{1i}\e}\cdot\cos(q\e\sqrt{\lambda\e})-\l
u\e\r_{S_{2i}\e}\bigg)\cdot w_1(x_i\e)+\\+\suml_{i\in
\I}{(d\e)^{N-1}\omega\sqrt{\lambda\e}\over\sin(q\e\sqrt{\lambda\e})}\bigg(\l
u\e\r_{S_{2i}\e}\cdot\cos(q\e\sqrt{\lambda\e})-\l
u\e\r_{S_{1i}\e}\bigg)\cdot w_2(x_i\e)+ \delta\e
\end{multline}
and the remainder
\begin{multline}\label{remainder}
\delta\e=\suml_{k=1,2}\intl_{\Omega_k\e}\ds\suml_{i\in
\I}\bigg[\nabla\bigg(\big(w_k(x_i\e)-w_k(x)\big)\phi_i\e(x)\bigg)\cdot\nabla
u\e(x)-\\-\lambda\e \big(w_k(x_i\e)-w_k(x)\big)\phi_i\e(x)
u\e(x)\bigg] dx
\end{multline}
is vanishingly small as $\eps\to 0$ (since
$|w(x_i\e)-w\e(x)|<Cd\e$ for $x\in\mathrm{supp}(\phi_i\e)$\ ):
\begin{gather*}
|\delta\e|\leq C\|u\e\|^2_{H^1(\Omega)}\suml_{i\in
\I}\left|\mathrm{supp}(\phi_i\e) \right|\underset{\eps\to 0}\to 0.
\end{gather*}

We rewrite (\ref{plugg1}) in the form
\begin{multline}\label{plugg1+}
0=\suml_{k=1,2}\intl_{\Omega\e}\bigg(\nabla (\Pi_k\e
u\e)\cdot\nabla w_k -\lambda\e (\Pi_k\e u\e )
w_k\bigg)dx+\\+{(d\e)^{N-1}\omega\over\eps^N}\cdot{\sqrt{\lambda\e}\over
\sin(q\e\sqrt{\lambda\e})} \intl_{\Omega}\bigg(\widehat{S}_1\e
u\e\cdot \cos(q\e\sqrt{\lambda\e})-\widehat{S}_2\e
u\e\bigg)\widehat{X}\e w_1
dx+\\+{(d\e)^{N-1}\omega\over\eps^N}\cdot{\sqrt{\lambda\e}\over
\sin(q\e\sqrt{\lambda\e})} \intl_{\Omega}\bigg(\widehat{S}_2\e
u\e\cdot \cos(q\e\sqrt{\lambda\e})-\widehat{S}_1\e
u\e\bigg)\widehat{X}\e w_2 dx+\delta\e,
\end{multline}
where the operators $\widehat{X}\e,\ \widehat{S}_k\e\ (k=1,2)$ are
defined by the formulae (\ref{oper1}).

It is obvious that for any $w\in C^1(\Omega)$\ $\widehat{X}\e
w\underset{\eps\to 0}\to w$ strongly in $L_2(\Omega)$. Moreover
since $p, q>0$ then $\D=\infty$ and therefore due to Corollary
\ref{cor1} (see (\ref{in_cor1_2})) $\widehat{S}_k\e
u\e\underset{\eps\to 0}\to u_k$ strongly in $L_2(\Omega)$.
Therefore passing to the limit (as $\eps\to\ 0$) in
(\ref{plugg1+}) we conclude that
\begin{multline}\label{conclud1}
0=\suml_{k=1,2}\intl_{\Omega}\bigg(\nabla u_k\cdot\nabla w_k
-\lambda_0 u_k w_k\bigg)dx+{p\omega\over q}{\sqrt{\lambda_0}\over
\sin(q\sqrt{\lambda_0})}
\intl_{\Omega}\bigg(u_1\cdot\cos(q\sqrt{\lambda_0})-u_2\bigg)w_1
dx+\\+{p\omega\over q}{\sqrt{\lambda_0}\over
\sin(q\sqrt{\lambda_0})}
\intl_{\Omega}\bigg(u_2\cdot\cos(q\sqrt{\lambda_0})-u_1\bigg)w_2dx,\quad
\forall w_1,w_2\in C_0^\infty(\Omega).
\end{multline}
It is easy to see that (\ref{conclud1}) implies that
$\mathrm{A}(\lambda_0)\mathrm{U}=0$. The fulfilment  (\ref{ah}) is
proved.

\textbf{Step 2.} Let us prove the fulfilment of the condition
(\ref{bh}) of the Hausdorff convergence.

Firstly suppose that $\lambda_0\in\sigma(\mathrm{A}(\lambda))$. We
have to prove that there exists $\lambda\e\in\sigma(-\Delta\e)$
such that $\lambda\e\underset{\eps\to 0}\to\lambda_0$.

Proving this indirectly we assume the opposite. Then the
subsequence (still denoted by $\eps$) and a positive number
$\delta$ exist such that
\begin{gather}\label{dist1}
\minl_{\lambda\e\in
\sigma(-\Delta\e)}|\lambda_0-\lambda\e|>\delta.
\end{gather}

Since $\lambda_0$ belongs to the spectrum of $\mathrm{A}(\lambda)$
there exists
$\mathrm{F}=\left(\begin{matrix}f_1\\f_2\end{matrix}\right)\in
[L_2(\Omega)]^2$ such that
\begin{gather}\label{notinIm}
\mathrm{F}\not\in \mathrm{Im}\mathrm{A}(\lambda_0).
\end{gather}

Let us consider the following problem on $M\e$:
\begin{gather}\label{B-problem1}
-\Delta\e u - \lambda_0 u=f\e.
\end{gather}
In view of (\ref{dist1}) this problem has the unique solution
$u\e(\tilde x)\in H_0^1(M\e)$ for an arbitrary $f\e\in L_2(M\e)$.
We set
\begin{gather*}
f\e(\tilde x)=\begin{cases}f_k(x),&\tilde x\in (x,k),\
k=1,2,\\
0,&\tilde x\in G_i\e,\ i\in\I.
\end{cases}
\end{gather*}

One has
\begin{gather*}
\|u\e\|_{L_2(M\e)}\leq {\|f\e\|_{L_2(M\e)}\over\delta}\leq C_1,\\
\|\nabla\e u\e\|^2_{L_2(M\e)}\leq |\lambda_0|\cdot
\|u\e\|_{L_2(M\e)}^2+\left|(f\e,u\e)_{L_2(M\e)}\right|\leq C_2.
\end{gather*}
Therefore $\|\Pi_k\e u\e\|_{H_1(\Omega)}\leq C$ ($k=1,2$) and by
the embedding theorem there exists a subsequence (still denoted by
$\eps$) such that $\Pi_k\e u\e\underset{\eps\to 0}\to u_k\in
H_0^1(\Omega)\quad (k=1,2)$.

For an arbitrary $w\e\in H^1_0(M\e)$ we have the following
equality:
\begin{gather}\label{int_ineq_f}
\intl_{M\e}\left(\nabla\e u\e\cdot\nabla\e w\e -\lambda\e u\e
w\e-f\e w\e\right) d\tilde x=0.
\end{gather}

Let us substitute into (\ref{int_ineq_f}) the function $w\e$
defined by the formula (\ref{we}) and pass to the limit in
(\ref{int_ineq_f}) as $\eps\to 0$. Similarly to "Step 1" we prove
that
$$\mathrm{A}(\lambda_0)\mathrm{U}=\mathrm{F},\quad
\mathrm{U}=\left(\begin{matrix}u_1\\u_2\end{matrix}\right).$$ Thus
we obtain a contradiction to (\ref{notinIm}).

In order to complite the verification of the fulfilment of
(\ref{bh}) we have to prove that for any $\lambda_0\in \cupl_{n\in
\mathbb{N}}\left\{{(\pi n)^2  q^{-2}}\right\}$ there exists
$\lambda\e\in\sigma(-\Delta\e)$ such that
$\lambda\e\underset{\eps\to 0}\to\lambda_0$. But this fact follows
directly from the structure of $\sigma(\mathrm{A}(\lambda_0))$
(see (\ref{struc1})). Indeed (\ref{struc1}) implies that for some
$n\in \mathbb{N}$ $\lambda_0=\liml_{m\to\infty}\lambda_{m}^n$,\
$\lambda_{m}^n\in\sigma(\mathrm{A}(\lambda))$, while we just prove
that $\minl_{\lambda\e\in
\sigma(-\Delta\e)}|\lambda_{m}^n-\lambda\e|\underset{\eps\to 0}\to
0$, hence $\minl_{\lambda\e\in
\sigma(-\Delta\e)}|\lambda_0-\lambda\e|\underset{\eps\to 0}\to 0$.

Therefore it remains to prove the statement (\ref{struc1}). We
will prove  it on the final third step.

\textbf{Step 3.} First of all let us note that on the Step 2 it
was proved that the spectrum  $\sigma(\mathrm{A}(\lambda))$ of
$\mathrm{A}(\lambda)$ belongs to $[0,\infty)$ (because each point
of $\sigma(\mathrm{A}(\lambda))$ is a limit of positive numbers
from $\sigma(-\Delta\e)$). Therefore now we are interested only in
the case $\lambda\geq 0$.

Let $\mathrm{U}=\left(\begin{matrix}u_1\\u_2\end{matrix}\right)\in
\left\{\mathrm{U}\in [H^2(\Omega)]^2,\
\mathrm{U}|_{\partial\Omega}=0\right\}$. We denote $u^{\pm}=u_1\pm
u_2$. Then it is easy to obtain that if
\begin{gather*}
\mathrm{A}(\lambda)\mathrm{U}=\mathrm{F},\quad
\mathrm{F}=\left(\begin{matrix}f_1\\f_2\end{matrix}\right)\in
[L_2(\Omega)]^2
\end{gather*}
then
\begin{gather*}
\begin{matrix}
-\Delta u^+-\left(\lambda+\ds{p\omega \over
q}\sqrt{\lambda}\tan\ds\left({q\sqrt{\lambda}\over
2}\right)\right)u^+=f^+,\\-\Delta u^--\left(\lambda-\ds{p\omega
\over q}\sqrt{\lambda}\cot\ds\left({q\sqrt{\lambda}\over
2}\right)\right)u^-=f^-,
\end{matrix}\quad f^{\pm}=f_1\pm f_2.
\end{gather*}
Thus $\sigma(\mathrm{A}(\lambda))$ coincides with the spectrum of
the pencil $\mathrm{\widetilde{A}}(\lambda)$ which is defined by
the operation
\begin{gather}\label{A_tilde}
\mathrm{\widetilde{A}}(\lambda)=\left(\begin{matrix}-\Delta
-\ds{p\omega\over q}\sqrt\lambda\tan \left({q\sqrt{\lambda}\over
2}\right)&0\\0&-\Delta +\ds{p\omega\over q}\sqrt\lambda\cot
\left({q\sqrt{\lambda}\over 2}\right)
\end{matrix}\right)-\lambda \mathrm{I}
\end{gather}
and by  the definitional domain
$\mathcal{D}(\mathrm{\widetilde{A}}(\lambda))=\mathcal{D}(\mathrm{{A}}(\lambda))$.
Obviously the spectrum of $\mathrm{\widetilde{A}}(\lambda)$
consists of such $\widehat{\lambda}\geq 0$ that solves at least
one of the following equations:
\begin{gather}
\label{tan} \widehat{\lambda}+{p\omega \over
q}\sqrt{\widehat{\lambda}}\tan\ds\left({q\sqrt{\widehat{\lambda}}\over
2}\right)=\mu\in\{\mu_m\}_{m\in
\mathbb{N}},\\\label{cot}\widehat{\lambda}- {p\omega \over
q}\sqrt{\widehat{\lambda}}\cot\ds\left({q\sqrt{\widehat{\lambda}}\over
2}\right)=\nu\in\{\mu_m\}_{m\in \mathbb{N}},
\end{gather}
where
$0<\mu_1\leq\mu_2\leq{\dots}\leq\mu_m\leq{\dots}\underset{m\to\infty}\to
\infty$ is the sequence of eigenvalues of the operator $-\Delta$
in $\Omega$ (with Dirichlet boundary conditions on
$\partial\Omega$).

Let $\mathcal{J}_n=\ds\left({(\pi (n-1))^2q^{-2}},{(\pi
n)^2q^{-2}}\right)$\ where $n$ be odd. Then it is easy to obtain
that if $m>M_n$, where $M_n$ is sufficiently large number
depending on $n$, then in $\mathcal{J}_n$ the equation (\ref{tan})
with the right-hand-side $\mu=\mu_m$ has the unique root
${\lambda}_{m}^{n,\tan}$ and moreover
${\lambda}_{m}^{n,\tan}\underset{m\to\infty}\to {(\pi n)^2q^{-2}
}$. The equation (\ref{cot}) also can have the roots
${\lambda}_{m}^{n,\cot}$ on the segment $\mathcal{J}_n$ but the
number of such roots is \textit{finite} (because on
$\mathcal{J}_n$ the function in the left-hand-side of (\ref{cot})
is bounded above). Note that possibly some $\widehat{\lambda}$
solve the equations (\ref{tan}) and (\ref{cot}) simultaneously (of
course in this case $\mu\not=\nu$).

Thus we have a countable set of point in $\mathcal{J}_n$ which are
the roots of one of the equations (\ref{tan}), (\ref{cot}) and
therefore this points are the eigenvalues of
$\mathrm{A}(\lambda)$. This set has only one accumulation point
$\ds{(\pi n)^2q^{-2}}$.

The same arguments are used if $n$ is even (in this case\
$\tan\left({q\sqrt{\lambda}\over 2}\right)$\ and\
$-\cot\left({q\sqrt{\lambda}\over 2}\right)$\ change places). Thus
the statement (\ref{struc1}) is proved that completes the proof of
Theorem \ref{th1}.\quad $\square$

\subsection{Proof of Theorem \ref{th2}}

The fulfilment of the condition (\ref{ah}) of the Hausdorff
convergence is proved similarly to that one in Theorem \ref{th1}.
Therefore we give only a sketch of the proof.

Let $\lambda\e\in\sigma(-\Delta\e)$, $\lambda\e\underset{\eps\to
0}\to\lambda_0$. If $\lambda_0\in \cupl_{n\in
\mathbb{N}}\left\{{(\pi n)^2q^{-2}}\right\}$ then (\ref{ah}) is
proved. So we consider the case $\lambda_0\notin \cupl_{n\in
\mathbb{N}}\left\{{(\pi n)^2q^{-2}}\right\}$.

Let $u\e$ be the eigenfunction that corresponds to $\lambda\e$ and
$\|u\e\|_{L_2(M\e)}=1$. Then there exists a subsequence (still
denoted by $\eps$) such that $\Pi_k\e u\e\underset{\eps\to 0}\to
u_k\in H_0^1(\Omega)$ $(k=1,2)$ strongly in $L_2(\Omega)$ and
weakly in $H^1(\Omega)$.

For an arbitrary $w\e\in H^1_0(M\e)$ the equality (\ref{int_ineq})
holds. We substitute into (\ref{int_ineq}) the function $w\e$ of
the form (\ref{we}). Using (\ref{in_cor1_3}) we pass to the limit
as $\eps\to 0$ in (\ref{int_ineq}) and obtain that
\begin{gather*}
\suml_{k=1,2}\intl_{\Omega}\left(\nabla u_k\nabla w_k-\lambda_0
u_kw_k\right)dx=0,\quad \forall w_1,w_2\in C_0^\infty(\Omega).
\end{gather*}
By Lemma \ref{lm3}
$\ds\suml_{i\in\I}\|u\e\|^2_{L_2(G_i\e)}\underset{\eps\to 0}\to 0$
and therefore
$\mathrm{U}=\left(\begin{matrix}u_1\\u_2\end{matrix}\right)\not=
0$. Thus $\lambda_0$ is the eigenvalue of the operator
$\mathrm{A}$ (\ref{A3}).

The fulfilment of the condition (\ref{bh}) in the case
$\lambda_0\notin\cupl_{n\in \mathbb{N}}\left\{(\pi n)^2
q^{-2}\right\}$ is proved completely similarly to that one in
Theorem \ref{th1}. Therefore it remains to verify the fulfilment
of (\ref{bh}) in the case $\lambda_0\in\cupl_{n\in
\mathbb{N}}\left\{(\pi n)^2 q^{-2}\right\}$.

Proving (\ref{bh}) indirectly we assume the opposite. Then the
subsequence (still denoted by $\eps$) and a positive number
$\delta$ exist such that (\ref{dist1}) holds.

Since $\lambda_0\in\cupl_{n\in \mathbb{N}}\left\{(\pi n)^2
q^{-2}\right\}$ then there exists $f\in L_2(0,q)$ such that the
problem
\begin{gather}\label{contr}
-u^{\prime\prime}-\lambda_0 u=f,\ x\in(0,q),\quad u(0)=u(q)=0
\end{gather}
has no solutions.

In order to simplify our calculations we suppose that $q\e\leq q$,
in the general case the proof needs some simple modifications.

For each $\eps$ we fix the number
$\mathbf{j}=\mathbf{j}(\eps)\in\I$ (we select this number
arbitrarily). We consider the problem (\ref{B-problem1}) on $M\e$
with $f\e\in L_2(M\e)$ defined by the formula
\begin{gather*}
f\e(\tilde x)=\begin{cases}(d\e)^{-{N-1\over 2}}f(z),&\tilde
x=(\phi,z)\in
G_{\mathbf{j}}\e,\\
0,&\text{otherwise}.
\end{cases}
\end{gather*}
In view of (\ref{dist1}) this problem has the unique solution
$u\e(\tilde x)\in H_0^1(M\e)$. Moreover since
$\|f\e\|^2_{L_2(M\e)}=\omega\|f\|^2_{L_2(0,q\e)}<C$ and by virtue
of (\ref{dist1}) the functions $u\e(\tilde x)$ are bounded in
$H^1(M\e)$ uniformly in $\eps$.

On $G_{\mathbf{j}}\e$ we represent $u\e$  in the form
$u\e(\phi,z)=\overline{u\e}(z)+v\e(\phi,z)$, where
$\overline{u\e}(z)=\l u\e\r_{S_{\mathbf{j}}[z]}$\ (recall that
$S_{\mathbf{j}}\e[z]=\left\{\tilde x=(\phi,\tau)\in
G_\mathbf{j}\e:\ \tau=z\right\}$). Notice that due to the Poincare
inequality
\begin{gather}
\label{v} \|v\e\|^2_{L_2(G_\mathbf{j}\e)}\leq
C(d^\eps)^2\|\nabla\e
u\e\|^2_{L_2(G_\mathbf{j}\e)}\underset{\eps\to 0}\to 0
\end{gather}

We introduce the operator $\Pi\e: H^1(G_{\mathbf{j}}\e)\to
H^1(0,q)$ by the formula
\begin{gather*}
\Pi\e u\e(z)=(d\e)^{{N-1\over
2}}\cdot\begin{cases}\overline{u\e}(z),&z\in [0,q\e],\\
\overline{u\e}(q\e),&z\in [q\e,q].\end{cases}
\end{gather*}
Due to the Cauchy inequality we have:
\begin{gather}\label{est1}
\|\Pi\e u\e\|^2_{L_2(0,q)}\leq
{1\over\omega}\|u\e\|^2_{L_2(G_{\mathbf{j}}\e)}+(q-q\e)\big| \Pi\e
u\e(q\e)\big|^2,\quad \|(\Pi\e u\e)_z^\prime\|^2_{L_2(0,q)}\leq
{1\over\omega}\|\nabla\e u\e\|^2_{L_2(G_{\mathbf{j}}\e)}.
\end{gather}
Furthermore using fundamental theorem of calculus it is easy to
obtain that
\begin{gather}\label{est2}
\big|\Pi\e u\e(q\e)\big|^2\leq 2\left[(q\e)^{-1}\|\Pi\e
u\e\|_{L_2(0,q\e)}^2+q\e\intl_0^{q\e} \left|{\partial \Pi\e
u\e\over\partial z}\right|^2dz\right].
\end{gather}
It follows from (\ref{est1})-(\ref{est2}) that the functions
$\Pi\e u\e$ are bounded in $H^1(0,q)$ uniformly in $\eps$.
Therefore there exists a subsequence (still denoted by $\eps$)
such that
\begin{gather*}
\Pi\e u\e\underset{\eps\to 0}\to u\in H^1(0,q)\text{ strongly in
}L_2(0,q)\text{ and weakly in }H^1(0,q).
\end{gather*}

We have the estimate
\begin{gather}\label{est3}
\l u\r^2_{S_{ki}\e}\leq C\left({\|\nabla\e
u\e\|^2_{L_2(B_{ki}\e)}\over \D\e}+{\|u\e\|^2_{L_2(B_{ki}\e)}\over
\eps^{N}}\right),\ k=1,2,\ i\in \I,
\end{gather}
which is proved similarly to (\ref{in_lm1_1}). Using (\ref{est3})
and the trace theorem we obtain
\begin{gather*}
|u(0)|^2=\liml_{\eps\to 0}\bigg((d\e)^{N-1}\l
u\r^2_{S_{1\mathbf{j}}\e}\bigg)\leq C\liml_{\eps\to
0}\left({(d\e)^{N-1}\|\nabla\e
u\e\|^2_{L_2(B_{k\mathbf{j}}\e)}\over
\D\e}+{(d\e)^{N-1}\|u\e\|^2_{L_2(B_{k\mathbf{j}}\e)}\over
\eps^{N}}\right)= 0.
\end{gather*}
Similarly $u(q)=0$. Thus $u\in H_0^1(0,q)$.

Let $w\in C^\infty(0,q)$ be an arbitrary function  such that
$\supp(w)\subset [\delta,q-\delta]$, where $\delta=\delta(w)$ is
some positive number. Since $q\e\underset{\eps\to 0}\to q$ then
for sufficiently small $\eps$ $\supp(w)\subset
[\delta,q\e-\delta/2]$. We define $w\e\in H^1_0(M\e)$ by the
formula:
\begin{gather*}
w\e(\tilde x)=\begin{cases}(d\e)^{-{N-1\over 2}}w(z),&\tilde
x=(\phi,z)\in
G_{\mathbf{j}}\e,\\
0,&\text{otherwise}.
\end{cases}
\end{gather*}

Then we have
\begin{gather}\notag
0=\liml_{\eps\to 0}\intl_{M\e}\left(\nabla\e u\e\cdot\nabla\e
w\e-\lambda_0 u\e w\e -f\e w\e\right)d\tilde x= \liml_{\eps\to
0}\intl_{G_{\mathbf{j}}\e}\left(\nabla\e
\overline{u\e}\cdot\nabla\e w\e-\lambda_0 \overline{u\e} w\e -f\e
w\e\right)d\tilde
 x+\delta(\eps)=\\\label{int_ineq1+}=\omega\intl_0^q \left({d\Pi\e
u\e(z)\over dz} {d w(z)\over dz}-\lambda_0 \Pi\e u\e(z) w(z)-
f(z)w(z)\right)dz+\delta(\eps),
\end{gather}
where the reminder\
$\delta\e=\ds\intl_{G_{\mathbf{j}}\e}\left(-v\e\Delta\e
w\e-\lambda_0 v\e w\e\right)d\tilde x$\ tends to zero as $\eps\to
0$ by virtue of (\ref{v}) and the definition of $w\e$. Passing to
the limit as $\eps\to 0$ in (\ref{int_ineq1+}) we obtain
\begin{gather*}
\intl_0^q \left({d u(z)\over dz} {d w(z)\over dz}-\lambda_0 u(z)
w(z)- f(z)w(z)\right)dz=0
\end{gather*}
for any function $w\in C^\infty(0,q)$ such that $\supp(w)\subset
[\delta,q-\delta]$, where $\delta=\delta(w)>0$ is some positive
number. Since the set of such functions is dense in $H^1_0(0,q)$
we conclude that $u$ is the solution to (\ref{contr}). We obtain a
contradiction.

Thus the fulfilment of (\ref{bh}) is completely verified. Theorem
\ref{th2} is proved.

\subsection{Proof of Theorem \ref{th3}}

We restrict ourselves to the proof of fulfilment of condition
(\ref{ah}). The condition (\ref{bh}) is proved using the same idea
as in Theorem \ref{th1}.

So let $\lambda\e\underset{\eps\to 0}\to\lambda_0$, $u\e$ be the
corresponding eigenfunction such that $\|u\e\|_{L_2(M\e)}=1$. Then
there exists a subsequence still denoted by $\eps$ such that
$\Pi_k\e u\e\underset{\eps\to 0}\to u_k\in H^1_0(\Omega)\ (k=1,2)$
strongly in $L_2(\Omega)$ and weakly in $H^1(\Omega)$. Due to
(\ref{in_cor1_4}) and (\ref{in_cor1_2}): $u_1=u_2$ and
$\liml_{\eps\to 0} \|S_k\e u\e-u\|^2_{L_2(\Omega)}= 0$ (here we
denote $u=u_k$).

For an arbitrary $w\e\in H^1_0(M\e)$ we have the equality
(\ref{int_ineq}). Let us introduce the following test function
$w\e$:
\begin{gather*}
w\e(\tilde x)=\begin{cases}w(x)+\ds\suml_{i\in
\I}\bigg(w(x_i\e)-w(x)\bigg)\cdot\phi\left({|x-x_i\e|\over
d\e}\right),&\tilde x=(x,k)\in \Omega_k\e,\ k=1,2,\\
w(x_i\e),&\tilde x=(\phi,z)\in G_i\e,\ i\in\I.
\end{cases}
\end{gather*}
Here $w\in C_0^\infty(\Omega)$ is an arbitrary function,
$\phi(r):[0,\infty)\to \mathbb{R}$ is a smooth positive function
equal to $1$ as $r\leq 1$ and equal to $0$ as $r\geq 2$.

Substituting this function into (\ref{int_ineq}) we obtain that
\begin{gather}\label{int_ineq1}
\suml_{k=1,2}\intl_{\Omega\e}\bigg(\nabla(\Pi_k\e u\e)\cdot\nabla
w-\lambda\e(\Pi_k\e u\e) w\bigg)dx-\lambda\e\omega\suml_{i\in\I}(d
\e)^{N-1}w(x_i\e)\intl_0^{q\e} \l
u\e\r_{S_{i}\e[z]}dz+\delta(\eps),
\end{gather}
where the remainder $\delta(\eps)$ has the form (\ref{remainder})
and tends to zero as $\eps\to 0$. Also we have
\begin{gather*}
\lambda\e\omega\suml_{i\in\I}(d\e)^{N-1}w(x_i\e)\intl_0^{q\e} \l
u\e\r_{S_{i}\e[z]}dz=\lambda\e\omega{(d\e)^{N-1}q\e\over
\eps^N}\intl_\Omega \widehat{X}\e w\cdot \widehat{S}_k\e u\e
dx+\delta_1(\eps),\quad k=1\vee k=2,
\end{gather*}
where the reminder $\delta_1(\eps)$ tends to zero by virtue of the
inequity (\ref{in_lm1_2}).

Passing to the limit in (\ref{int_ineq1}) we conclude that
\begin{gather}\label{int_ineq2}
\intl_\Omega\bigg(2\nabla u\cdot\nabla w-\lambda_0(2+p\omega)
uw\bigg)dx=0,\quad\forall w\in C_0^\infty(\Omega).
\end{gather}
By virtue of (\ref{in_cor1_4})\ \ $1=\liml_{\eps\to
0}\|u\e\|^2_{L_2(M\e)}=(2+p\omega)\|u\|^2_{L_2(\Omega)}$.
Therefore $u\not=0$.

Thus (\ref{int_ineq2}) implies that $\lambda_0$ is the eigenvalue
of the operator defined by the operation $-\left(1+{1\over
2}{p\omega}\right)^{-1}\Delta$ and Dirichlet boundary conditions.
Theorem \ref{th3} is proved.

\subsection{Proof of Theorem \ref{th4}}
As in the previous theorem we restrict ourselves to the proof of
fulfilment of the condition (\ref{ah}).

Let $\lambda\e\underset{\eps\to 0}\to\lambda_0$, $u\e$ be the
corresponding eigenfunction such that $\|u\e\|^2_{L_2(M\e)}=1$.
Then there exists a subsequence (still denoted by $\eps$) such
that $\Pi_k\e u\e\underset{\eps\to 0}\to u_k\in H^1_0(\Omega),\
k=1,2$ (strongly in $L_2(\Omega)$ and weakly in $H^1(\Omega)$).

By Lemma \ref{lm2}:\ $1=\liml_{\eps\to 0}\|u\e\|^2_{L_2(M\e)}=
\|u_1\|^2_{L_2(\Omega)}+\|u_2\|^2_{L_2(\Omega)}$.

For an arbitrary $w\e\in H^1_0(M\e)$ the equality (\ref{int_ineq})
holds.

We start from the case $\D=\infty,\ r<\infty$. Let us consider the
following test function $w\e$:
\begin{gather*}
w\e(\tilde x)=\begin{cases}w_k(x)+\ds\suml_{i\in
\I}\bigg(w_k(x_i\e)-w_k(x)\bigg)\cdot\phi\left({|x-x_i\e|\over
d\e}\right),&\tilde x=(x,k)\in \Omega_k\e,\ k=1,2,\\
v_i\e(z),&\tilde x=(\phi,z)\in G_i\e,\ i\in\I.
\end{cases}
\end{gather*}
Here $w_k\in C_0^\infty(\Omega)$ ($k=1,2$) are arbitrary
functions, $\phi(r):[0,\infty)\to \mathbb{R}$ is a smooth positive
function equal to $1$ as $r\leq 1$ and equal to $0$ as $r\geq 2$,
$v_i\e(z)$ is defined by the formula
$$v_i\e(z)=z\cdot\ds{w_2(x_i\e)-w_1(x_i\e)\over q\e}+w_1(x_i\e)$$

Substituting this function into (\ref{int_ineq}) we obtain that
\begin{gather}\notag
0=\suml_{k=1,2}\intl_{\Omega\e}\bigg(\nabla(\Pi_k\e
u\e)\cdot\nabla w_k-\lambda\e(\Pi_k\e u\e)
w\bigg)dx+\suml_{i\in\I}(d\e)^{N-1}\omega\left({\partial
v_i\e(q\e)\over\partial z}\l u\e \r_{S_{2i}\e}-{\partial
v_i\e(0)\over\partial z}\l u\e
\r_{S_{1i}\e}\right)+\\\notag+\left(\delta\e+\delta_1\e\right)=
\suml_{k=1,2}\intl_{\Omega\e}\bigg(\nabla(\Pi_k\e u\e)\cdot\nabla
w_k-\lambda\e(\Pi_k\e u\e)
w_k\bigg)dx+\\\label{int_ineq3}+{(d\e)^{N-1}\omega\over \eps^N
q\e}\intl_\Omega\left(\widehat{S}_1\e u\e-{\widehat{S}}_2\e
u\e\right)\left(\widehat{X}\e w_1-\widehat{X}\e
w_2\right)dx+\left(\delta\e+\delta_1\e\right),
\end{gather}
where the remainder $\delta(\eps)$ has the form (\ref{remainder})
and tends to zero as $\eps\to 0$, the remainder $\delta_1(\eps)$
has the form
\begin{gather}\label{remainder1}
\delta_1\e=-\lambda\e\suml_{i\in \I}\intl_{G_i\e}u\e w\e d\tilde x
\end{gather}
and tends to zero in view of Lemma \ref{lm2}.

Passing to the limit as $\eps\to 0$ in (\ref{int_ineq3}) we obtain
the following equality:
\begin{gather}\notag
0=\suml_{k=1,2}\intl_{\Omega}\big(\nabla u_k\cdot\nabla
w_k-\lambda_0 u_k w_k\big)dx+r\omega
\intl_{\Omega}(u_1-u_2)(w_2-w_2)dx,\quad\forall w_1,w_2\in
C_0^\infty(\Omega).
\end{gather}
This equality implies that $\lambda_0$ is the eigenvalue of the
operator $\mathrm{A}$ (\ref{A2}).

Now we consider the case $\D<\infty$. For simplicity we restrict
ourselves to the case $N>2$, in the case $N=2$ the proof is
completely similar.

We substitute into the equality (\ref{int_ineq}) the following
test function:
\begin{gather*}
w\e(\tilde x)=\begin{cases}w_k(x)+\ds\suml_{i\in
\I}\left(\bigg(w_k(x_i\e)-w_k(x)\bigg)\cdot\phi\left({|x-x_i\e|\over
d\e}\right)+\right.\\\qquad\qquad\left.+ \bigg(v_i\e(\tilde
x)-w_k(x_i\e)\bigg)\cdot\phi\left(4\cdot\ds{|x-x_i\e|\over
\eps}\right)\right),
&\tilde x=(x,k)\in \Omega_k\e,\ k=1,2,\\
v_i\e(\tilde x),&\tilde x=(\phi,z)\in G_i\e,\ i\in\I.
\end{cases}
\end{gather*}
Here $w_k\in C_0^\infty(\Omega)$ ($k=1,2$) are arbitrary
functions, $\phi(r):[0,\infty)\to \mathbb{R}$ is a smooth positive
function equal to $1$ as $r\leq 1$ and equal to $0$ as $r\geq 2$,
$v_i\e(\tilde x)$ is the solution of the following problem:
\begin{gather*}
\begin{matrix}
-\Delta\e v_i\e=0,\ \tilde x\in G_i\e\cup B_{1i}\e\cup
B_{2i}\e,\qquad v_i\e=w_k(x_i\e),\ \tilde x\in C_{ki}\e
\end{matrix}
\end{gather*}
(the sets $B_{ki}\e,\ C_{ki}\e$ are defined at the beginning of
Section \ref{sec2}). It is easy to calculate $v_i\e$:
\begin{gather}\label{v_formula}
v_i\e=\begin{cases}{a_{ki}\e\cdot
|x-x_i\e|^{2-N}}+b_{ki}\e,&\tilde x=(x,k)\in\Omega_k\e,\
k=1,2,\\A_i\e z+B_i\e,&\tilde x=(\phi,z)\in G_i\e,\
i\in\I,\end{cases}
\end{gather}
where
$a_{1i}\e=-a_{2i}\e=\ds{(d\e)^{N-2}(w_2(x_i\e)-w_1(x_i\e))\over
2+\ds{q\e\over
d\e}(N-2)-2\left({2d\e\over\eps}\right)^{N-2}}=A_i\e\ds{(d\e)^{N-1}\over
N-2}$, $b_{ki}\e=w_k(x_i\e)-a_{ki}\e \ds\left({\eps\over
2}\right)^{2-N}$, $B_i\e=a_{1i}\e (d\e)^{2-N}+b_{1i}\e$.

Integrating by parts in (\ref{int_ineq}) we obtain:
\begin{multline}\label{int_ineq4}
0=-\suml_{k=1,2}\intl_{\Omega\e}\bigg(\Pi_k\e u\e\cdot\Delta
w+\lambda\e(\Pi_k\e u\e)
w\bigg)dx-\\-\ds\suml_{k=1,2}\suml_{i\in\I} \intl_{\Omega\e}
\Delta\left(\bigg(w_k(x_i\e)-w_k(x)\bigg)\cdot\phi\left({|x-x_i\e|\over
d\e}\right)\right)u\e dx-\\-\ds\suml_{k=1,2}\suml_{i\in\I}
\intl_{\Omega\e}
\Delta\left(\bigg(v_i\e(x)-w_k(x_i\e)\bigg)\cdot\phi\left(4\cdot{|x-x_i\e|\over
\eps}\right)\right)u\e dx-\lambda\e\suml_{i\in\I}\intl_{G_i\e}u\e
w\e d\tilde x.
\end{multline}
In view of Lemma \ref{lm2} the last term in (\ref{int_ineq4})
tends to zero as $\eps\to 0$. The second term (we denote it
$\delta\e$) can be rewritten in the form
\begin{gather}\notag
\delta\e=\ds\suml_{k=1,2}\suml_{i\in\I} \left(\intl_\Omega
\nabla\left(\bigg(w_k(x_i\e)-w_k(x)\bigg)\cdot\phi\left({|x-x_i\e|\over
d\e}\right)\right)\cdot \nabla\left( \Pi_k\e u\e\right)
dx-\intl_{D_i\e}\Delta w_k\cdot \Pi_k\e u\e dx\right)
\end{gather}
and converges to zero as $\eps\to 0$ since
$\ds\suml_{i\in\I}\left|D_i\e\right|\underset{\eps\to 0}\to 0$ and
$\ds\suml_{i\in\I}\left|\supp (\phi_i\e)\right|\underset{\eps\to
0}\to 0$ (here $\phi_i\e=\phi\left({|x-x_i\e|/d\e}\right)$).

Finally we investigate the third term (we denote them
$\mathrm{I}\e$). One has:
\begin{gather*}
\mathrm{I}\e=-\ds\suml_{k=1,2}\suml_{i\in\I}\left(
\intl_{\Omega\e}
\Delta\left(\bigg(v_i\e(x)-w_k(x_i\e)\bigg)\cdot\phi\left({|x-x_i\e|\over
\eps}\right)\right)\left(u\e-\l u\e\r_{B_{ki}\e}\right) dx+ \l
u\e\r_{B_{ki}\e}\intl_{S_{ki}\e}\ds{\partial v_i\e\over\partial
\vec{n}}dx\right),
\end{gather*}
where ${\vec{n}}$ is the exterior normal to $S\e_{ki}$. Taking in
account (\ref{v_formula}) and Poincare inequality we obtain:
\begin{gather*}
\left|\ds\suml_{k=1,2}\suml_{i\in\I} \intl_{\Omega\e}
\Delta\left(\bigg(v_i\e(x)-w_k(x_i\e)\bigg)\cdot\phi\left({|x-x_i\e|\over
\eps}\right)\right)\left(u\e-\l u\e\r_{B_{ki}\e}\right)
dx\right|^2\leq\\\leq
C\eps^2\ds\left({\D\e\over\eps^N}\right)^2\suml_{k=1,2}\suml_{i\in\I}\|\nabla\e
u\e\|^2_{L_2(B_{ki}\e)}\underset{\eps\to 0}\to 0,
\end{gather*}
and
\begin{gather*}
-\ds\suml_{k=1,2}\suml_{i\in\I}\l
u\e\r_{B_{ki}\e}\intl_{S_{ki}\e}\ds{\partial v_i\e\over\partial
\vec{n}}dx=\suml_{i\in\I}{ (d\e)^{N-2}\omega (N-2)\left(\l
u\e\r_{B_{1i}\e}-\l u\e\r_{B_{2i}\e}\right)
\left(w_1(x_i\e)-w_2(x_i\e)\right)\over 2+\ds{q\e\over
d\e}(N-2)-2\left({2d\e\over\eps}\right)^{N-2}}=\\={
(d\e)^{N-2}\omega(N-2)\ds\intl_\Omega \left(\widehat{B}_1\e
u\e-\widehat{B}_2\e u\e\right)\left(\widehat{X}\e
w_1-\widehat{X}\e w_2\right)dx\over \eps^N\left(2+\ds{q\e\over
d\e}(N-2)-2\left({2d\e\over\eps}\right)^{N-2}\right)}\underset{\eps\to
0}\to \mathrm{V}\intl_\Omega (u_1-u_2)(w_1-w_2)dx.
\end{gather*}
where $\mathrm{V}$ is defined by (\ref{VVV}). We conclude that
\begin{gather*}
-\suml_{k=1,2}\intl_\Omega\left(u_k\Delta w_k+\lambda_0 u_k
w_k\right)dx+\mathrm{V}\intl_\Omega
(u_1-u_2)(w_1-w_2)dx,\qquad\forall w_1,w_2\in C_0^\infty(\Omega),
\end{gather*}
and thus $\lambda_0$ is the eigenvalue of the operator
$\mathrm{A}$ (\ref{A2}). Theorem \ref{th4} is proved.

\section{\label{sec4}Number-by-number convergence of eigenvalues and
convergence of eigenfunctions}

In the last section we study the convergence as $\eps\to 0$ of the
eigenvalue $\lambda_m\e$ for fix number $m\in \mathbb{N}$. Also we
describe the behavior of the eigenfunctions $u_m\e$.

We start from the case $q=0$. Let $\{\lambda_m\}_{m\in
\mathbb{N}}$ be the sequence of eigenvalues of homogenized
operator $\mathrm{A}$ that acts in $L_2(\Omega)$ and is defined by
the formula (\ref{A1}) if\ \ $(r=\infty\wedge \D=\infty)$ or acts
in $[L_2(\Omega)]^2$ and is defined by the formula (\ref{A2}) if\
\ $((r<\infty\wedge\D=\infty)\vee \D<\infty)$. The eigenvalues
$\lambda_m,\ m\in \mathbb{N}$ are renumbered in the increasing
order and are repeated according to their multiplicity. By
$N(\lambda_m)$ we denote the eigenspace that corresponds to
$\lambda_{m}$.

\begin{theorem}{}\label{th6} Let $q=0$. Then $\forall m\in
\mathbb{N}:\ \ \lambda_m\e\underset{\eps\to 0}\to\lambda_m.$
\end{theorem}

\noindent\textbf{Proof.} We present the proof, for example, in the
case $((r<\infty\wedge\D=\infty)\vee \D<\infty)$ (i.e. under the
conditions of Theorem \ref{th4}). For the case $(r=\infty\wedge
\D=\infty)$ theorem is proved in a similar way. The proof is based
on the following

\begin{lemma}\label{lm4}
Let $\widehat{\lambda}$ be the eigenvalue of homogenized operator
$\mathrm{A}$, let $\widehat{M}$ be the multiplicity of
$\widehat{\lambda}$. Suppose that for $j={m,{\dots},m+M-1}$\
$\liml_{\eps\to 0}\lambda_j\e=\widehat{\lambda}$ and
$\liml_{\eps\to 0}\lambda_{m-1}\e<\widehat{\lambda}<\liml_{\eps\to
0}\lambda_{m+M}\e$. Then $M=\widehat{M}$.
\end{lemma}

\noindent\textbf{Proof.} When proving Theorem \ref{th4} we show
that there exists a subsequence (still denoted by $\eps$) such
that
\begin{gather*}
\Pi_k u_j\e\underset{\eps\to 0}\to u_{kj}\e\in H^1_0(\Omega)\
(j=m,{\dots},m+M-1,\ k=1,2)\text{ strongly in }L_2(\Omega)\text{
and weakly in }H_0^1(\Omega)
\end{gather*}
and
$\mathrm{U}_j=\left(\begin{matrix}u_{1j}\\u_{2j}\end{matrix}\right)$
are the eigenfunctions of $\mathrm{A}$ which correspond to
$\widehat{\lambda}$. By Lemma \ref{lm2}
\begin{gather*}
\delta_{\a\b}=\liml_{\eps\to
0}(u_{\a}\e,u_{\b}\e)_{L_2(M\e)}=(\mathrm{U}_{\a},\mathrm{U}_{\b})_{[L_2(\Omega)]^2},\quad
\a,\b=m,{\dots},m+M-1.
\end{gather*}
So we have $M$ functions $\mathrm{U}_j$ ($j=m,{\dots},m+M-1$) that
belong to $N(\widehat{\lambda})$ and are orthonormal in
$[L_2(\Omega)]^2$. Hence $M\leq \widehat{M}$.

Now we prove that $M=\widehat{M}$. Assuming the opposite we
suppose that $M<\widehat{M}$. Let $\mathrm{H}$ be the subspace of
$N(\widehat{\lambda})$ generated by $\mathrm{U}_{j}$,
$j=m,{\dots},m+M-1$. By the assumption
$N(\widehat{\lambda})\ominus \mathrm{H}\not=\{0\}$.

Let
$\mathrm{F}=\left(\begin{matrix}f_{1}\\f_{2}\end{matrix}\right)\in
N(\widehat{\lambda})$. Then $\mathrm{F}\not\in
\mathrm{Im}(\mathrm{A}-\widehat{\lambda}\mathrm{I})$. We introduce
the function $f\e\in L_2(M\e)$:
\begin{gather*}
f\e(\tilde x)=\begin{cases}f_k(x),&\tilde x=(x,k)\in\Omega_k\e,\
k=1,2,\\0,&\tilde x\in \cupl_{i\in\I}G_i\e.\end{cases}
\end{gather*}
Let us consider the following problem
\begin{gather*}
-\Delta\e u-\widehat{\lambda} u=\widehat{f}\e,
\end{gather*}
where
$\widehat{f}\e=f\e-\ds\suml_{j=m}^{m+M-1}(f\e,u_j\e)_{L_2(M\e)}$.
For sufficiently small $\eps$\ \
$|\widehat{\lambda}-\lambda_i\e|\geq\delta>0$ if
$j\not=m,{\dots},m+M-1$. Therefore this problem has the unique
solution $u\e(\tilde x)\in H^1_0(\Omega)$ that is defined by the
formula
\begin{gather*}
u\e=\suml_{j\in \mathbb{N}:\
j\not=\overline{m,m+M-1}}{(\widehat{f}\e,u_j\e)_{L_2(M\e)}\over
\lambda_j\e-\widehat{\lambda}}u_j\e
\end{gather*}
moreover $$\|u\e\|_{L_2(M\e)}\leq
\delta^{-1}\|\widehat{f}\|_{L_2(M\e)}\leq C_1,\quad \|\nabla\e
u\e\|^2_{L_2(M\e)}\leq
\widehat{\lambda}\|u\e\|^2_{L_2(M\e)}+\left|(\widehat{f}\e,u\e)_{L_2(M\e)}\right|\leq
C_2.
$$
Hence the subsequence (still denoted by $\eps$) exists such that
$\Pi_k u\e\underset{\eps\to 0}\to u_{k}\in H^1_0(\Omega),\ k=1,2.$
In the same way as in Theorem \ref{th4} we conclude that
\begin{gather}\label{F-F}
(\mathrm{A}-\widehat{\lambda} \mathrm{I})
\mathrm{U}=\mathrm{F}-\suml_{j=\overline{m,m+M-1}}(\mathrm{F},\mathrm{U}_j)_{[L_2(\Omega)]^2},\quad
\mathrm{U}=\left(\begin{matrix}u_{1}\\u_{2}\end{matrix}\right).
\end{gather}
Now let us choose $\mathrm{F}$ from $N(\widehat{\lambda})\ominus
\mathrm{H}$. Then the right-hand-side in (\ref{F-F}) is equal to
$\mathrm{F}$ and therefore $\mathrm{F}\in
\mathrm{Im}(\mathrm{A}-\widehat{\lambda}\mathrm{I})$. We obtain a
contradiction. Lemma is proved.\qquad $\square$\medskip

It is easy to complete the proof of theorem. Let $\lambda_1$ has
the multiplicity $M_1$ (i.e.
$\lambda_1=\lambda_2={\dots}=\lambda_{M_1}<\lambda_{M_1+1}$). It
follows from the condition (\ref{bh}) of Hausdorff convergence
that $|\lambda\e_j|<C_j$ ($j=1,{\dots},M_1$). Let
$\eps'\subset\eps$ by an arbitrary subsequence such that
\begin{gather}\label{sub}
\exists\liml_{\eps'\to
0}\lambda_j^{\eps'}=\widehat{\lambda}_j,\quad j=1,{\dots},M_1.
\end{gather}
By the condition (\ref{ah}) of Hausdorff convergence
$\widehat{\lambda}_j\in\sigma(\mathrm{A})$. By the property
(\ref{bh}) $\widehat{\lambda}_1=~\lambda_1$. By Lemma \ref{lm4}
$\widehat{\lambda}_j=\lambda_j$ for $j=2,{\dots},M_1$. Since
$\eps'$ is an arbitrary subsequence for which (\ref{sub}) holds
then $\liml_{\eps\to 0}\lambda_j^{\eps}={\lambda}_j,\
j=1,{\dots},M_1$.

For the next $\lambda_j$ ($j>M_1$) the theorem is proved by
induction.\quad $\Box$

\begin{theorem}{}\label{th7}Let $q=0$. Let
$\lambda_{m-1}<\lambda_{m}=\lambda_{m+1}={\dots}=\lambda_{m+M-1}<\lambda_{M+m}$.

Then for any $w\in N(\lambda_{m})$ the linear combination
$\widehat{u}\e=\ds\suml_{j=m}^{m+M-1}\a_ju_j\e$ and the
subsequence $\eps'\subset\eps$ exist such that
\begin{gather}\label{vek_conv1}
\Pi^{\eps'}_k \widehat{u}^{\eps'}\underset{\eps'\to 0}\to w_k\
(k=1,2)\text{ strongly in }L_2(\Omega)\text{ and weakly in
 }H^1(\Omega),\end{gather}
where $w_1=w_2=w$ if\ \ $(r=\infty\wedge \D=\infty)$ or
$\left(\begin{matrix}w_1\\w_2\end{matrix}\right)=w$ if\ \
$((r<\infty\wedge\D=\infty)\vee \D<\infty)$.

Conversely for any linear combination of the form
$\widehat{u}\e=\ds\suml_{j=m}^{m+M-1}\a_ju_j\e$ there exist $w\in
N(\lambda_{m})$ and the subsequence $\eps'\subset\eps$ such that
(\ref{vek_conv1}) holds.

\end{theorem}

\noindent\textbf{Proof.} Let for distinctness
$((r<\infty\wedge\D=\infty)\vee \D<\infty)$ and the operator
$\mathrm{A}$ has the form (\ref{A2}) (for the case
$(r=\infty\wedge \D=\infty)$ the proof is similar). Then there
exists a subsequence $\eps'\subset\eps$ such that $\Pi_k
u_j\e\underset{\eps'\to 0}\to u_{kj}\e,\ j=m,{\dots},m+M-1,\
k=1,2$, the functions
$\mathrm{U}_j=\left(\begin{matrix}u_{1j}\\u_{2j}\end{matrix}\right)$
belong to $N(\lambda_m)$ and are orthonormal in $[L_2(\Omega)]^2$
(see the proof of Lemma \ref{lm4}). Then
$\left\{\mathrm{U}_j\right\}_{j=m}^{m+M-1}$ is the basis in
$N(\lambda_m)$ and therefore $w$ can be represented in the form
$w=\ds\suml_{j=m}^{m+M-1}\a_j\mathrm{U}_j$. We set
$\widehat{u}\e=\ds\suml_{j=m}^{m+M-1}\a_ju_j\e$. Obviously
$\widehat{u}\e$ satisfies (\ref{vek_conv1}).

Converse assertion actually is obtained within the proof of Lemma
\ref{lm4}.\quad $\square$ \medskip

Now let us consider the case $q>0$, $p>0$. Let
$0<\lambda_{1}^1\leq\lambda_{2}^1\leq{\dots}\leq\lambda_{m}^1\leq{\dots}\underset{m\to\infty}
\to \left\{\pi^2q^{-2}\right\}$ be the subsequence of eigenvalues
of operator pencil $\mathrm{A}(\lambda)$ (see Theorem \ref{th1})
that belong to the segment $(0,\pi^2 q^{-2})$. Here we write them
in increasing order and with account of their multiplicity.

\begin{theorem}{}\label{th8} Let $q>0$, $p>0$. Then $\forall m\in
\mathbb{N}:\ \lambda_m\e\underset{\eps\to 0}\to\lambda_{m}^1.$
\end{theorem}

\noindent\textbf{Proof.} The proof is directly follows from
Theorem \ref{th1} and from the following lemma which is an analog
of Lemma \ref{lm4}.

\begin{lemma}\label{lm5}
Let $\widehat{\lambda}$ be the eigenvalue of
$\mathrm{A}(\lambda)$, let $\widehat{M}$ be the multiplicity of
$\widehat{\lambda}$. Suppose that for $j={m,{\dots},m+M-1}$\
$\liml_{\eps\to 0}\lambda_j\e=\widehat{\lambda}$ and
$\liml_{\eps\to 0}\lambda_{m-1}\e<\widehat{\lambda}<\liml_{\eps\to
0}\lambda_{m+M}\e$. Then $M=\widehat{M}$.
\end{lemma}

\noindent\textbf{Proof.} When proving Theorem \ref{th1} we show
that there exist a subsequence (still denoted by $\eps$) such that
\begin{gather*}
\Pi_k u_j\e\underset{\eps\to 0}\to u_{kj}\e\in H^1_0(\Omega)\
(j=m,{\dots},m+M-1,\ k=1,2)\text{ strongly in }L_2(\Omega)\text{
and weakly in }H_0^1(\Omega)
\end{gather*}
and
$\mathrm{U}_j=\left(\begin{matrix}u_{1j}\\u_{2j}\end{matrix}\right)$
are the eigenfunctions of the pencil $\mathrm{A}(\lambda)$ that
correspond to $\widehat{\lambda}$.

Using Lemma \ref{lm3} and the parallelogram identity we obtain for
$\forall\alpha,\beta\in\{m,{\dots},m+M-1\}$:
\begin{multline}\label{scal}
\delta_{\a\b}=\liml_{\eps\to
0}(u_\a\e,v_\b\e)_{L_2(M\e)}=\\=\left(\mathrm{U}_\a,\mathrm{U}_\b\right)_{[L_2(\Omega)]^2}+
k_1\left(\mathrm{U}_\a,\mathrm{U}_\b\right)_{[L_2(\Omega)]^2}+
k_2\left((u_{1\a},u_{2\b})_{L_2(\Omega)}+(u_{2\a},u_{1\b})_{L_2(\Omega)}\right).
\end{multline}

Let us recall that also $\widehat{\lambda}$ is the eigenvalue of
the pencil $\mathrm{\widetilde{A}}(\lambda)$ (\ref{A_tilde}), it
solves one of the equations (\ref{tan}),(\ref{cot}) (possibly it
solves both equations). By $\widetilde{N}(\widehat{\lambda})$ we
denote the corresponding eigenspace, obviously $\mathrm{dim}
\widetilde{N}(\widehat{\lambda})=\mathrm{\dim}
N(\widehat{\lambda})$. We denote $u_j^{\pm}=u_{1j}\pm u_{2j}$.
Then the functions
$\mathrm{\widetilde{U}}_j=\left(\begin{matrix}u_{j}^+\\u_{j}^-\end{matrix}\right)$
belong to $\widetilde{N}(\widehat{\lambda})$.

One has $$\forall \a,\b\in m,{\dots},m+M-1:\quad
(u_\a^+,u_\b^-)_{L_2(\Omega)}=0.$$ Indeed if $\widehat{\lambda}$
solves only one of the equations (\ref{tan}) and (\ref{cot}) then
either $u_\a^+=0$ or $u_\b^-=0$ while if $\widehat{\lambda}$
solves both the equations (\ref{tan}), (\ref{cot}) then $u_\a^+$
and $u_\b^-$ are the eigenfunction of the operator $-\Delta$
corresponding to some (nonequal\ !) eigenvalues $\mu^+$ and
$\nu^-$.

Then it is easy to rewrite (\ref{scal}) in the form
\begin{gather*}
\delta_{\a\b}=\liml_{\eps\to
0}(u_\a\e,u_\b\e)_{L_2(M\e)}=\rho^+(u_\a^+,u_\b^+)+\rho^-(u_\a^-,u_\b^-),\quad
\rho^\pm={1\over 2}+{1\over 2}(k_1\pm k_2).
\end{gather*}
Remark that $\rho^\pm>\ds{1\over 2}$ since $k_1\pm
k_2=\ds{\left(1\mp\cos\left(q\sqrt{\widehat{\lambda}}\right)\right)
\left(q\sqrt{\widehat{\lambda}}\pm\sin\left(q\sqrt{\widehat{\lambda}}\right)\right)\over
2q\sqrt{\widehat{\lambda}}\sin^2\left(q\sqrt{\widehat{\lambda}}\right)}>0$.
Therefore
$\delta_{\a\b}=\left(\mathrm{\widetilde{U}}_\a,\mathrm{\widetilde{U}}_\b\right)_{\mathrm{\widetilde{H}}}$,
where $\mathrm{\widetilde{H}}=L_2(\Omega,\rho^+ dx)\oplus
L_2(\Omega,\rho^- dx)$. Therefore the functions
$\mathrm{\widetilde{U}}_j$ ($j=m,{\dots},m+M-1$) are linearly
independent and thus $M\leq\dim
\widetilde{N}(\widehat{\lambda})=\widehat{M}$.

The proof that $M=\widehat{M}$ is completely similar to the proof
of this equality in Lemma \ref{lm4}.\quad $\square$

\begin{theorem}{}\label{th9}Let $q>0$, $p>0$. Let
$\lambda_{m-1}^1<\lambda_{m}^1=\lambda_{m+1}^1={\dots}=\lambda_{m+M-1}^1<\lambda_{M+m}^1$.

Then for any $w=\left(\begin{matrix}w_1\\w_2\end{matrix}\right)\in
N(\lambda_{m}^1)$\ the linear combination
$\widehat{u}\e=\ds\suml_{j=m}^{m+M-1}\a_ju_j\e$  and the
subsequence $\eps'\subset\eps$ exist such that (\ref{vek_conv1})
holds.

Conversely for any linear combination of the form
$\widehat{u}\e=\ds\suml_{j=m}^{m+M-1}\a_ju_j\e$ there exist
$w=\left(\begin{matrix}w_1\\w_2\end{matrix}\right)\in
N(\lambda_{m}^1)$ and the subsequence $\eps'\subset\eps$ such that
(\ref{vek_conv1}) holds.

\end{theorem}

The proof is similar to the proof of Theorem \ref{th7}. \medskip

And finally we consider the case $q>0$, $p=0$. Here we restrict
ourselves to the investigation of the number-by-number convergence
of eigenvalues. By $\{\lambda_m\}_{m\in \mathbb{N}}$ we denote the
sequence of eigenvalues of the operator $\mathrm{A}$ that acts in
$[L_2(\Omega)]^2$ and is defined by (\ref{A3}). The eigenvalues
$\{\lambda_m\}_{m\in \mathbb{N}}$ are renumbered in the increasing
order and are repeated according to their multiplicity. We denote:
$$\mathcal{M}=\maxl_{\lambda_m< \pi^2 q^{-2}} m$$

\begin{theorem}{}\label{th10} Let $q>0$, $p=0$. Then
$\lambda_m\e\underset{\eps\to 0}\to\lambda_{m}$ if $m\leq
\mathcal{M}$ and $ \lambda_m\e\underset{\eps\to 0}\to \pi^2
q^{-2}$ otherwise.
\end{theorem}

\noindent\textbf{Proof.} Using Lemma \ref{lm3} (see
(\ref{in_lm3_2})) one can easily proof that for any
$\widehat{\lambda}\in \{\lambda_m\}_{m\in \mathbb{N}}$ such that
$\widehat{\lambda}\notin \{(\pi n)^2 q^{-2}\}_{n\in \mathbb{N}}$
the assertion of Lemma \ref{lm4} holds true. Then in the same way
as in the proof of Theorem \ref{th6} we conclude that
$\lambda_m\e\underset{\eps\to 0}\to\lambda_{m}$ if $m\leq
\mathcal{M}$ and $\lambda\e_{\mathcal{M}+1}\underset{\eps\to 0}\to
\pi^2 q^{-2}$.

It remains to proof that $\lambda_m\e\underset{\eps\to 0}\to \pi^2
q^{-2}$ if $m>\mathcal{M}+1$. We prove this by induction. Suppose
that $\lambda_m\e\underset{\eps\to 0}\to \pi^2 q^{-2}$ if
$\mathcal{M}+1\leq m\leq \mathcal{M}+\mu$, $\mu>0$. Then we have
to prove that $\lambda_{\mathcal{M}+\mu+1}\e\underset{\eps\to
0}\to \pi^2 q^{-2}$.

By $\mathbf{j}=\mathbf{j}(\eps)$ we denote such multiindex
$\mathbf{j}=\mathbf{j}(\eps)\in\I$ that
$$
\suml_{\mathcal{M}+1}^{\mathcal{M}+\mu}\left\|u_m\e\right\|^2_{L_2(G_\mathbf{j}\e)}\leq
\suml_{\mathcal{M}+1}^{\mathcal{M}+\mu}\left\|u_m\e\right\|^2_{L_2(G_i\e)}\text{
for }\forall i\in \I$$ It is easy to see that
\begin{gather}\label{sum}
\suml_{\mathcal{M}+1}^{\mathcal{M}+\mu}\left\|u_m\e\right\|^2_{L_2(G_\mathbf{j}\e)}\leq
C\eps^N
\end{gather}

Let us introduce the following function ${v}\e\in H^1_0(M\e)$:
\begin{gather*}
v\e=\begin{cases}\left[{1\over
2}q\e\omega(d\e)^{N-1}\right]^{-1/2}\cdot\sin\left(\pi (q\e)^{-1}
z\right),& \widetilde x=(\phi,z)\in
G_{\mathbf{j}}\e,\\
0,&\text{otherwise}.
\end{cases}
\end{gather*}
We have:
\begin{gather}\label{v1}
\|v\e\|^2_{L_2(M\e)}=1,\quad \|\nabla\e v\e\|^2_{L_2(M\e)}=\pi^2 (q\e)^{-2},\\
\label{v2} \ (v\e,u_m\e)_{L_2(M\e)}\underset{\eps\to 0}\to 0,\quad
m=1,{\dots},\mathcal{M}+\mu.
\end{gather}
The statement (\ref{v2}) follows from (\ref{in_lm3_2}) for
$m={\overline{1,\mathcal{M}}}$ and from (\ref{sum}) for
$m=\overline{\mathcal{M}+1,\mathcal{M}+\mu}$.

We denote
$$
\overline{v}\e=v\e-\suml_{m=1}^{\mathcal{M}+\mu}u_m\e(v\e,u_m\e)_{L_2(M\e)}
$$
Since $(\overline{v}\e,u_m\e)_{L_2(M\e)}=0$ for
$m=1,{\dots},\mathcal{M}+\mu$ then by Courant minimax principle
\begin{gather}\label{v4}
\lambda_{\mathcal{M}+\mu+1}\e\leq {\|\nabla\e
\overline{v}\e\|^2_{L_2(M\e)}\over \|
\overline{v}\e\|^2_{L_2(M\e)}}
\end{gather}
It follows from (\ref{v1})-(\ref{v4}) that
$\limsup\limits_{\eps\to 0}\lambda_{\mathcal{M}+\mu+1}\e\leq \pi^2
q^{-2}$.

On the other hand $\pi^2 q^{-2}=\lim\limits_{\eps\to
0}\lambda_{\mathcal{M}+\mu}\e\leq\liminf\limits_{\eps\to
0}\lambda_{\mathcal{M}+\mu+1}\e$. Thus $\lim\limits_{\eps\to
0}\lambda_{\mathcal{M}+\mu+1}\e=\pi^2 q^{-2}$.

 Theorem is proved.\quad $\square$
\medskip

\section*{Acknowledgments}

The author is grateful to Prof. E.Ya.Khruslov for the attention he
paid to this work. Also I would like to thank Prof. T.A.Melnyk for
the fruitful discussion, especially in the case $q>0$, $p=0$. The
work is partially supported by the joint French-Ukrainian project
"PICS 2009-2011. Mathematical Physics: Methods and Applications".

\end{document}